\documentclass[3p,times]{elsarticle}
\usepackage{ecrc}
\usepackage{amsmath}
\DeclareMathOperator\arctanh{arctanh}
\usepackage{graphics,graphicx,epsfig}
\usepackage{subfig}
\usepackage{psfrag}
\usepackage{color}
\usepackage{float}
\usepackage{color}
\usepackage{verbatim}
\usepackage{ulem}

\volume{00}

\firstpage{1}

\journalname{Commun Nonlinear Sci Numer Simulat}

\runauth{Alexander G. Korotkov, Tatiana A. Levanova, Michael A. Zaks, Grigory V. Osipov}

\jid{procs}

\jnltitlelogo{Commun Nonlinear Sci Numer Simulat}

\usepackage{amssymb}
\usepackage{amsmath,amsthm}
\usepackage{amsfonts}
\usepackage{graphicx}

\begin{document}

\begin{frontmatter}



\dochead{}

\title{Dynamics in a phase model of half-center oscillator:
two neurons with excitatory coupling}

\author[1]{Alexander G. Korotkov}
\author[1]{Tatiana A. Levanova}
\author[2]{Michael A. Zaks}
\author[1]{Grigory V. Osipov}
\address[1]{Control Theory Department, Institute of Information Technologies, Mathematics and Mechanics, Lobachevsky University, Gagarin ave. 23, Nizhny Novgorod, 603950, Russia}
\address[2]{Institute of  Physics, Humboldt University of Berlin Newtonstr. 15, Berlin, D-12489, Germany}

\begin{abstract}

A minimalistic model of the half-center oscillator is proposed. Within it, we consider dynamics of two excitable neurons interacting by means of the excitatory coupling. In the parameter space of the model, we identify the regions of dynamics, characteristic for central pattern generators: respectively, in-phase, anti-phase synchronous oscillations and quiescence, and study various bifurcation transitions between all these states. Suggested model can serve as a building block of specific complex central pattern generators for studies of rhythmic activity and information processing in animals and humans.

\end{abstract}

\begin{keyword}
half-center oscillator\sep central pattern generator \sep  theta-neuron\sep in-phase spiking\sep anti-phase spiking\sep bifurcations


\end{keyword}

\end{frontmatter}

% \linenumbers

\section{Introduction}

Central pattern generators (CPGs) are circuits in self-contained integrative nervous systems, able to generate and control basic repetitive patterns of coordinated motor behaviour without sensory feedback or peripheral input. They are responsible for such vital rhythmic motor behaviours as heartbeat, respiratory functions and locomotion \cite{Selverston1985} - \cite{katz2007invertebrate}. One of the best-known case studies in this field is of  locomotion in vertebrates: several decades of evidence (see e.g. \cite{mackay2002central}) support the hypothesis that walking, flying, and swimming are largely governed by a small network of spinal neurons in all vertebrate species, from lampreys to humans. Recent evidence suggests that plasticity changes of some CPG elements may contribute to the development of specific pathophysiological conditions associated with impaired locomotion or spontaneous locomotor-like movements \cite{guertin2013central}. Despite the relevance of the topic and substantial progress in the field, including proposed mechanisms of pattern generation  \cite{matsuoka1987mechanisms}-\cite{pusuluri2020computational}, genesis of the motor patterns is still not fully understood \cite{selverston2000reliable}.

One of the most widespread approaches in the numerical modelling of CPGs (as well as of other neuronal networks) uses the Hodgkin-Huxley equations \cite{izhikevich2007dynamical} or different kinds of their reductions, such as the FitzHugh-Nagumo equations \cite{pusuluri2020computational}, delivering detailed description of CPG. 

Since reproduction of temporal patterns, not the dynamics of an individual neuron, plays a crucial role \cite{sakurai2011different} in the paradigm of CPG, one may use reduction to phase equations in order to lower the computational complexity. The patterns of motor activity find expression in robust evolution of phase differences between the network elements, therefore it looks reasonable to adopt a phase oscillator as a model of an individual neuron. This approach goes back to the early modelling of animal locomotor CPG, where coupled systems of ODEs were reduced to phase models \cite{cohen1982nature}--\cite{buono2001models}.

Our goal is a model of CPG based on simple neuron-like units, able, on the one hand, to emulate a number of CPG dynamical patterns observed in experiments and reproduced in biologically plausible models \cite{wojcik2014key}-\cite{jalil2013toward}, and amenable, on the other hand, to analytical studies. 

Biological experiments witness that most CPGs have some kind of a universal constituent known as a half-center oscillator (HCO) \cite{hill2003half}. To account for the generation of rhythmic pattern, Brown \cite{brown1911intrinsic} first proposed the concept of HCO, in which two mutually inhibitory coupled neurons burst in anti-phase. HCO can consist of endogenously bursting neurons, intrinsically tonic spiking ones or even quiescent neurons that start to generate alternating activity when coupled. As shown in numerous theoretical studies \cite{wang1992alternating}-\cite{terman2008reducing}, formation of anti-phase bursting rhythm is tightly connected to slow time scale dynamics, associated with the slow membrane currents. Simple HCO can contribute to more complex modular CPG networks, such as swimming CPG of \textit{Melibe leonina} and \textit{Dendronotus iris} \cite{alaccam2015making}.

To understand better the dynamical principles underlying the behaviours of larger networks, we introduce a simple model of HCO based on two coupled units. Individual element in this case is an active rotator described by the Adler equation:
\begin{equation} \label{Adler_eq}
\mathop{\phi}\limits^\cdot = \gamma - \sin \phi,
\end{equation}
where $\phi$ corresponds to the phase of the element and $\gamma$ is a control parameter. 
\begin{figure}[H]
	\begin{center}
		\includegraphics[width = 0.35\linewidth]{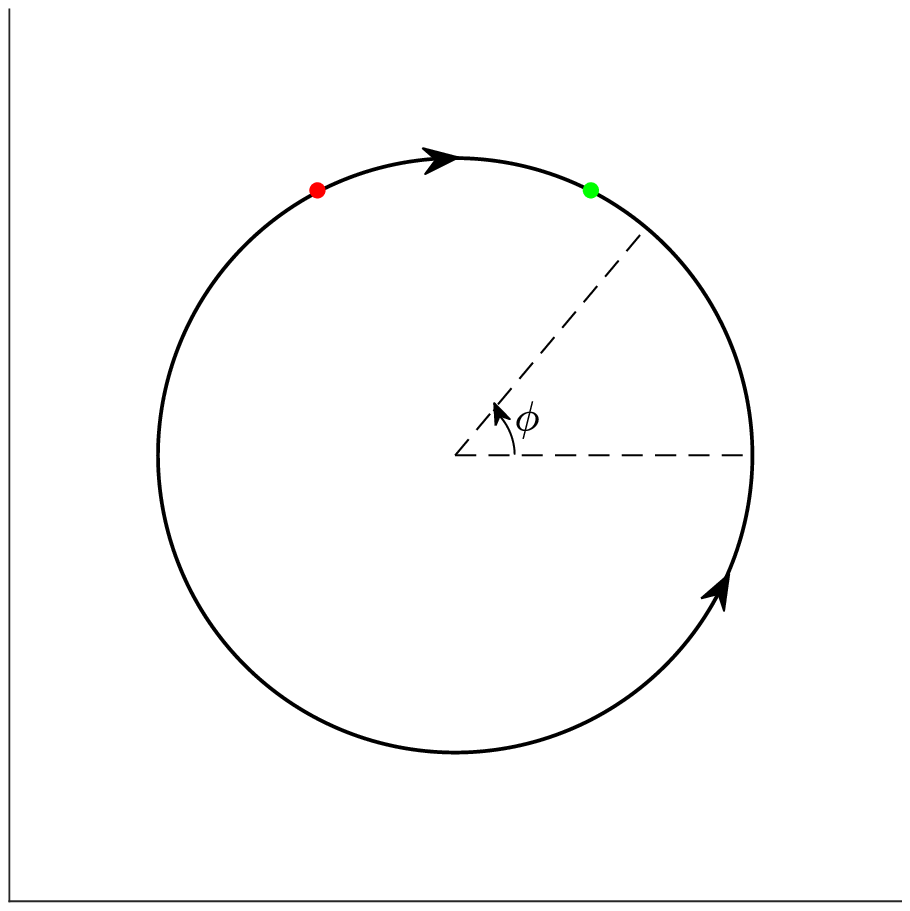}
		\includegraphics[width = 0.35\linewidth]{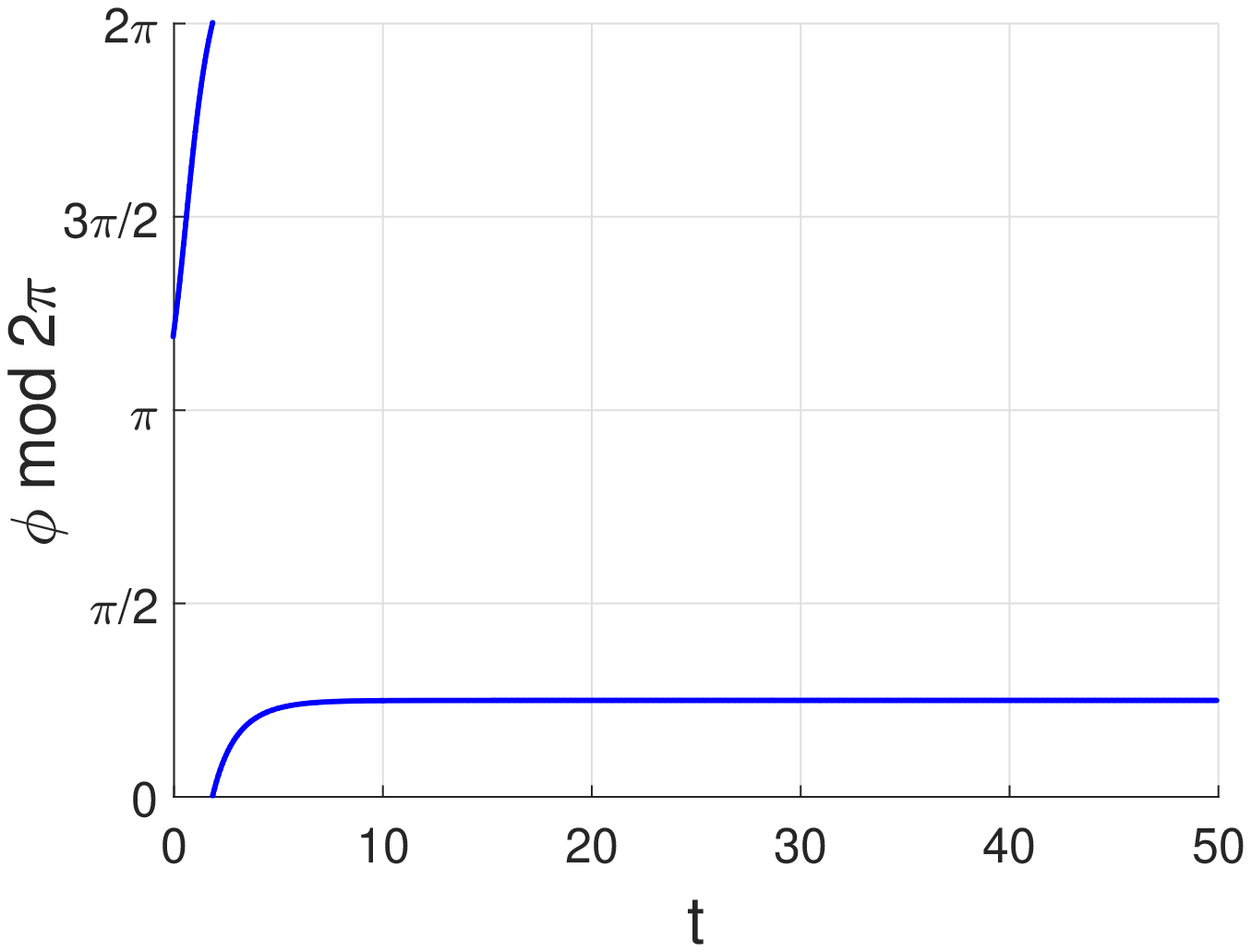}\\
		\includegraphics[width = 0.35\linewidth]{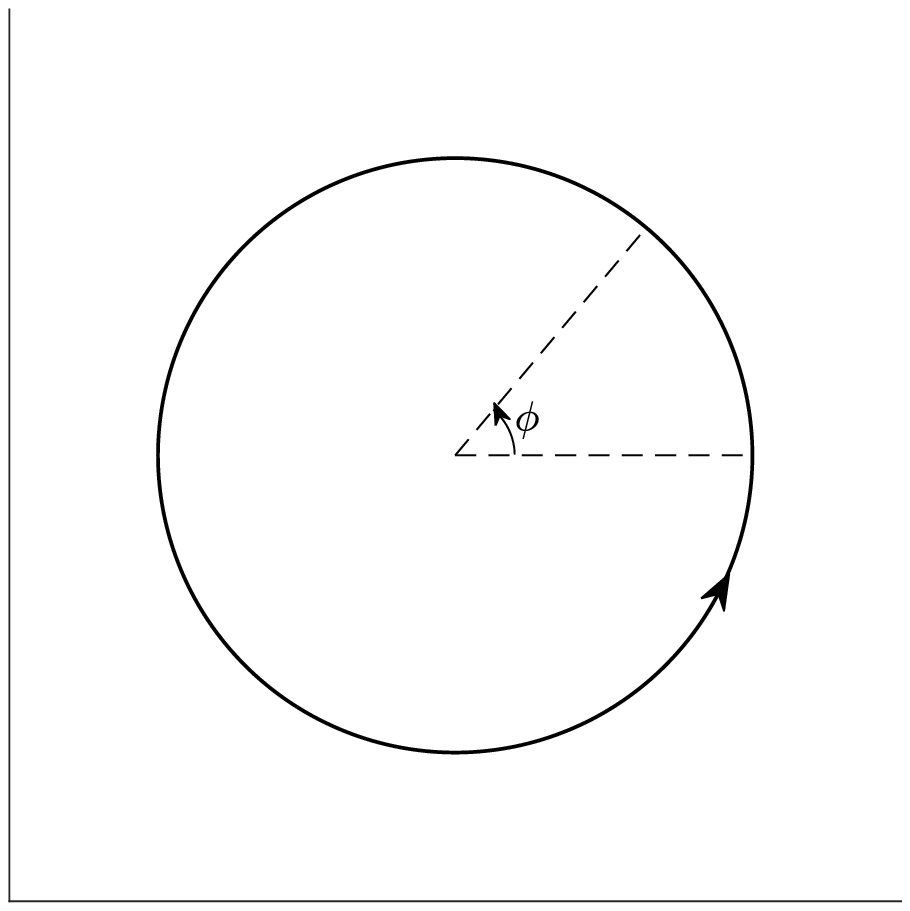}
		\includegraphics[width = 0.35\linewidth]{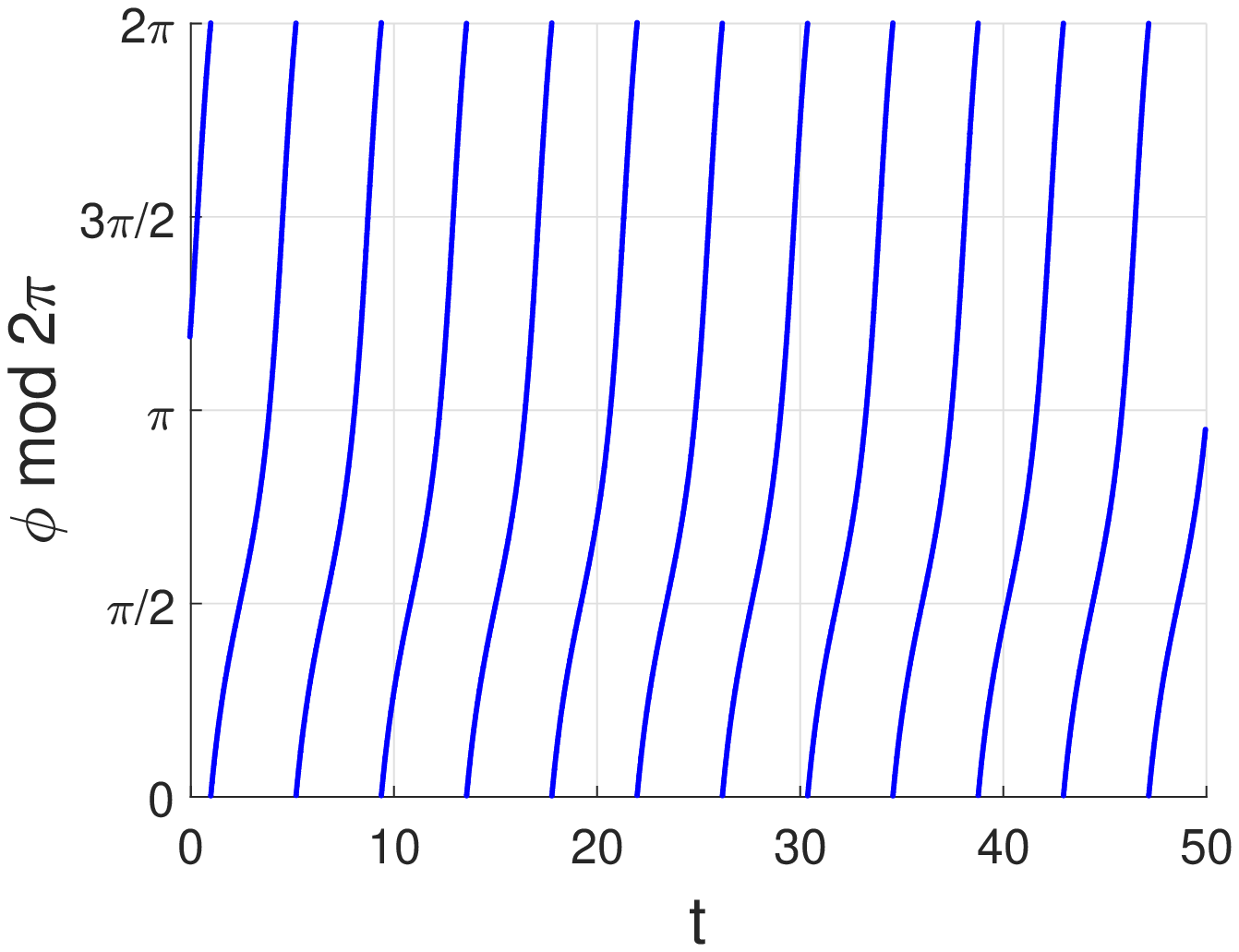}
	\end{center}
	\caption{ Phase space (circle) and time series of a single neuron-like element, described by Eq.\eqref{Adler_eq}. Top row:  Excitable dynamics at $\gamma<1$. Phase point is attracted to the state of rest (green dot on the circle on the left upper panel), which corresponds to the constant value of the phase $\phi$ of the element (right upper panel). %In this case Eq.\eqref{Adler_eq} describes excitable neuron. 
	Bottom row: Oscillatory dynamics  at $\gamma>1$. 
Phase of the single element evolves continuously in time (left bottom panel), and the element generates spikes (right bottom panel).  		\label{two_states}  }
\end{figure}

This model, introduced in \cite{adler1973study}, is evidently similar to the classical theta-neuron equation \cite{ermentrout1986parabolic}. In dependence on $\gamma$, Eq.\eqref{Adler_eq} reproduces excitable behaviour ($\gamma<1$, see upper panels in Fig.~\ref{two_states}) or self-oscillatory behaviour ($\gamma>1$, lower panels in Fig.~\ref{two_states}). Below we consider the first case.

In the present study our point is to understand, by studying symmetries and bifurcations, the basic types of dynamics typical for CPG and to gain more insights in the fundamental principles of HCO functioning that allow CPG to be ultimately flexible and multifunctional \cite{rubin2012explicit}-\cite{briggman2008multifunctional}.

The paper is organized as follows. First, we propose a simple phenomenological model of HCO and describe the way we have constructed it. Further, we discuss general properties of the introduced model. After that we focus on the main types of neuron-like activity typical for biological HCO. Our study concerns the properties of these states, as well as the bifurcation transitions between them. 
In conclusion, we summarize  our findings, discuss the directions of future studies.

\section{The simple model of HCO and its basic properties}

As a simple model of HCO we propose the motif of two identical excitable neurons, mutually interacting via the excitatory coupling. Mathematically it is described by a system of two differential equations:
\begin{equation} \label{ensemble}
\begin{cases}
\mathop{\phi_1}\limits^\cdot = \gamma - \sin \phi_1 + d\cdot I(\phi_2)\\
\mathop{\phi_2}\limits^\cdot = \gamma - \sin \phi_2 + d\cdot I(\phi_1)
\end{cases}.
\end{equation}
Here, the parameter $d$ ($d > 0$) regulates the strength of symmetric excitatory couplings $I(\phi)$.

\begin{figure}[H]
	\centering
	(a)\includegraphics[width = 0.45\linewidth]{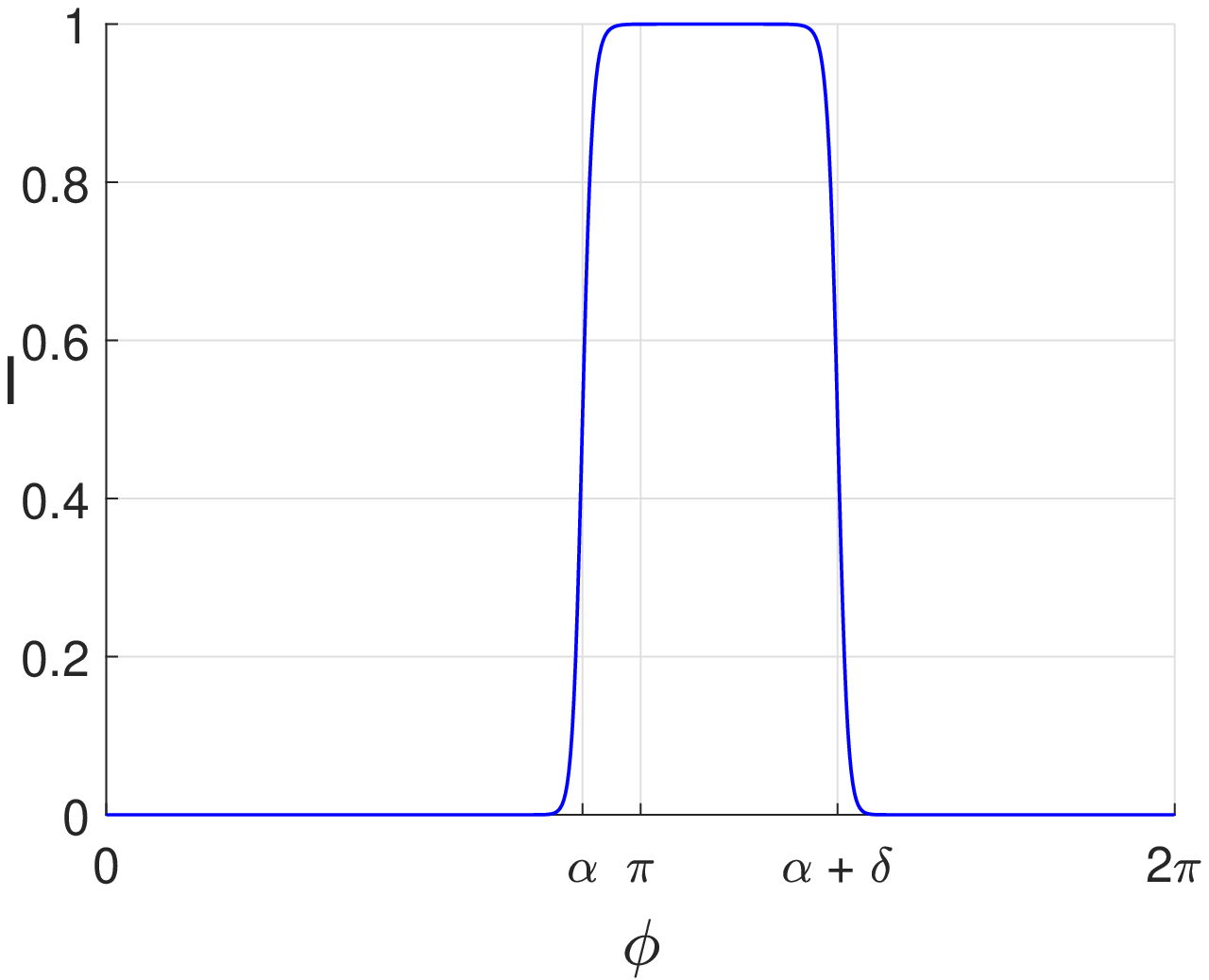}
	(b)\includegraphics[width = 0.45\linewidth]{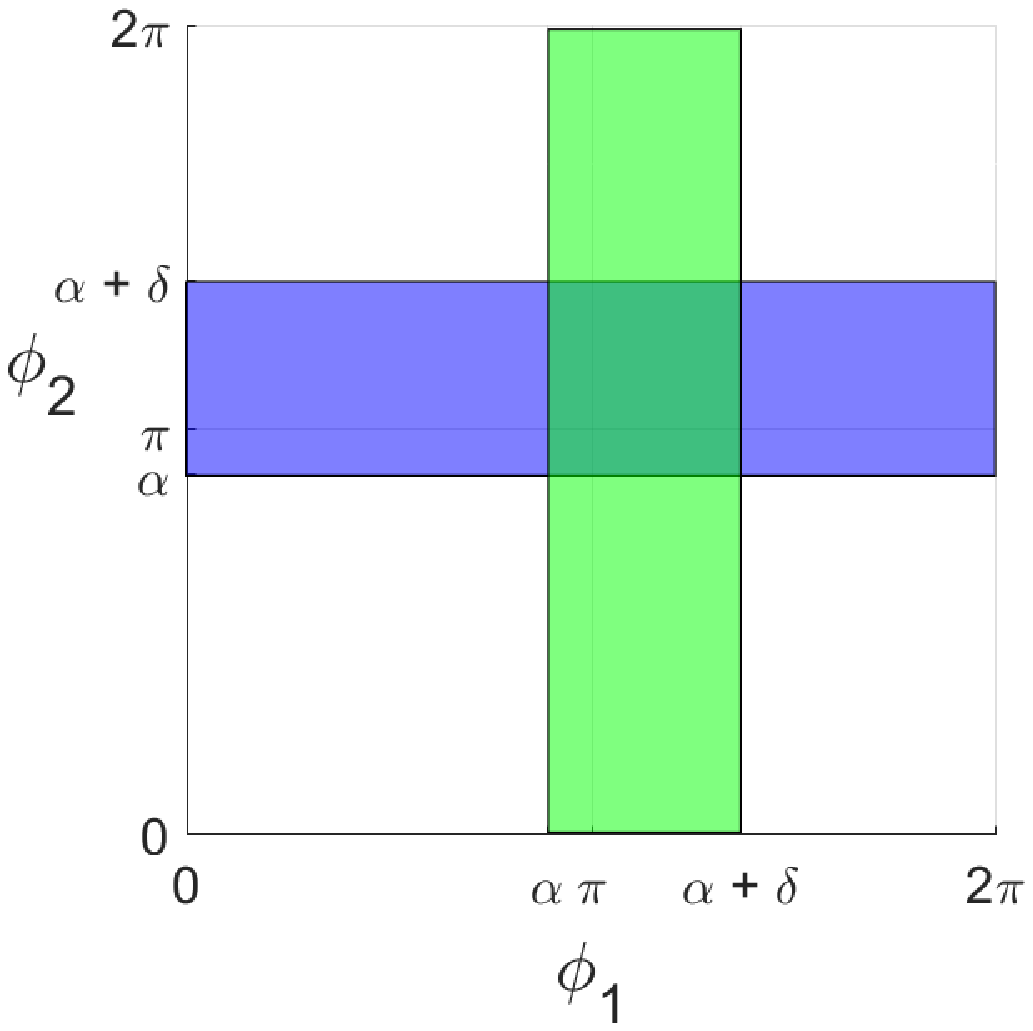}
	\caption{(a) Coupling function $I(\phi)$  (b) Regions of activation on the phase torus. Green coloured region: the first element activates the second one. Blue coloured region: the second element activates the first one. Mutual activation of both elements takes place in the overlap of two regions. Parameter values: $k=50$, $\alpha=\frac{7\pi}{8}$, $\delta=\frac{\pi}{2}$.}
	\label{coupling_function}
\end{figure}

In accordance to the biological principles \cite{destexhe1994efficient}, we model excitatory coupling by the function
\begin{equation} \label{coupling_func}
I(\phi) = \frac{1}{1 + e^{k\,\big(\cos(\delta/2) - \cos(\phi - \alpha - \delta/2)\big)}}.
\end{equation}
Coupling of this form, first introduced in \cite{korotkov2019dynamics}, and tested in subsequent studies~\cite{korotkov2018chaotic,korotkov2019effects}, simulates the transmission of a signal from the presynaptic element to the postsynaptic one. 
The coupling function \eqref{coupling_func} takes into account the basic principles of chemical synaptic coupling: (i) presence/absence of the activity in the postsynaptic element depends on the activity level in the presynaptic one; (ii) all interactions between neuron cells are inertial due to the fact that the transfer of neurotransmitter is not instantaneous.
When the phase $\phi$ of the active presynaptic element reaches the value $\alpha$, the current is applied to the postsynaptic element. Duration of the impact of this stimulus is defined by the difference $\delta$. The parameter $k$ regulates the steepness of transitions between the open and practically closed states of a synapse: the larger the value of $k$, the sharper are the transitions. Dependence of the coupling function $I(\phi)$ on the phase $\phi$ of the presynaptic element is sketched in Fig. \ref{coupling_function}(a). The diagram in Fig. \ref{coupling_function}(b) shows the regions of the joint phase space, where the elements are mutually activated. 

The phase space of Eq.~\eqref{ensemble} is a two-dimensional torus.
Altogether the system \eqref{ensemble} with the coupling \eqref{coupling_func} is governed by five parameters: 
$\gamma, d, k,\alpha,\delta$. Of these, we fix below the values $\gamma = 0.7$ (putting thereby the individual elements into the excitable states) and  $k = 50$ (this %large value 
ensures the sharp profile of $I(\phi)$).

The remaining active parameters $\alpha$ and $\delta$ are responsible for inertia and duration effects, respectively; by adjusting them, we can simulate synapses with different neurotransmitters.
Formally, the period of the coupling function with respect to the parameter $\delta$ is $4\pi $. In fact, $\delta$ takes values from the interval $[0, 2\pi)$, since the activation range is the segment $[\alpha, \alpha + \delta]$, that is, at $\delta = 2\pi$ both elements always activate each other.

As already mentioned, we focus both on various types of neuron-like activity, like the in-phase and anti-phase spiking patterns, and on bifurcation scenarios behind the onset and destruction of these patterns in the simple model  \eqref{ensemble} of the HCO. Below, the term \textit{in-phase limit cycle} refers to a limit cycle in which the phases of both elements coincide: $\phi_1(t) = \phi_2(t)$. Further, \textit{anti-phase limit cycle} denotes a limit cycle with some period $T$ in which the phases are shifted with regards to each other of by half-period: $\phi_1(t) = \phi_2(t + \frac{T}{2})$. 

Let us briefly discuss the basic features of the system \eqref{ensemble}, utilizable for further analysis.
We start with properties that hold regardless of the (continuous) function $I(\phi)$.

\textit{Property 1.} 
Since the system \eqref{ensemble} is invariant under a permutation of variables $\phi_{1,2}$, the phase portrait is symmetric with respect to the invariant diagonal $\phi_1 = \phi_2$.

\textit{Property 2.} Suppose that an anti-phase cycle exists in the phase space of the system \eqref{ensemble}. 
Then, for each of its points $(\phi_1^*, \phi_2^*)$, the cycle also contains the symmetrical counterpart 
$(\phi_2^*, \phi_1^*)$, shifted in time by the half-period of this cycle.

\textit{Property 3.} Two or more anti-phase limit cycles cannot coexist in the phase space of the system.

We start the proof of this property with a remark that an anti-phase cycle, due  to Property 2, cannot be entirely confined  either to the triangle $0<\phi_1<\phi_2<2\pi$ or to the symmetric triangle $0<\phi_2<\phi_1<2\pi$. Hence, the phase curve of the cycle should intersect the axes $\phi_1=0$ and $\phi_2=0$.

Assume that there are two anti-phase limit cycles. Let the first one include a point with coordinates $(0,a)$, where $0<a<2\pi$. Then (Property 2) it also contains a point with coordinates $(a,0)$, which on the 2-torus is identified with a point $(a,2\pi)$. Let the second anti-phase cycle pass through the points with coordinates $(0,b)$ and $(b,2\pi)$ ($0<b<2\pi$), and let $b$ exceed $a$. 
Two continuous curves crossing the triangle $0<\phi_1<\phi_2<2\pi$, so that the first of them passes through the points with coordinates $(0,a)$ and $(a,2 \pi)$, whereas the second contains points  $(0,b)$ and $(b,2\pi)$,
 are obliged to intersect. This invalidates the assumption on the existence of more than one anti-phase cycle.

\textit{Property 4.} The system \eqref{ensemble} has two types of equilibria: the equilibria of the first type lie on the line $\phi_1 = \phi_2$, the equilibria of the second type lie elsewhere and, due to the symmetry, appear in pairs with coordinates of the form $(a, b)$ and $(b, a)$.

Let us prove that the existence of a pair of equilibria of the second type implies the existence of a equilibrium of the first type. The coordinates of the latter $\phi_1 = \phi_2 = \phi$ fulfill the equation
\begin{equation} \label{equilibrium_sym}
\gamma - \sin \phi + d\cdot I(\phi) = 0.
\end{equation}
Similarly, equilibria of the second type can be recovered from the system
\begin{equation} \label{equilibriums}
\begin{cases}
\gamma - \sin \phi_1 + d\cdot I(\phi_2) = 0\\
\gamma - \sin \phi_2 + d\cdot I(\phi_1) = 0
\end{cases}.
\end{equation}
Suppose that a pair of equilibria of the second type $(a, b)$ and $(b, a)$ exists. 
%Therefore, one can investigate the following system of equations
Their coordinates $a,b$ are solutions of
\begin{equation*}
\begin{cases}
\gamma - \sin a + d\cdot I(b) = 0\\
\gamma - \sin b + d\cdot I(a) = 0
\end{cases},
\end{equation*}
whence follows $\gamma - \sin a + d\cdot I(a) = -(\gamma -\sin b + d\cdot I(b))$, i.e. function $F(\phi) = \gamma - \sin \phi + d\cdot I(\phi)$ takes values of different signs (or zeros) at $\phi = a$ and $\phi = b$. Then, by virtue of continuity, there exists $\xi$ ($a \leq \xi \leq b$) such that $\gamma - \sin \xi + d\cdot I(\xi) = 0$, i.e. $\xi$ satisfies \eqref{equilibrium_sym}. Thus, the existence of a pair of equilibrium states of the second type implies the existence of a equilibrium state of the first type.

Further properties concern the specific coupling function \eqref{coupling_func}.

\textit{Property 5.} The system \eqref{ensemble} is invariant under the transformation $\phi_i \longrightarrow \pi - \phi_i$, $t \longrightarrow -t$, $\alpha \longrightarrow \pi - \alpha - \delta$. 
It follows that the bifurcation diagram in the parameter space ($\alpha$, $\delta$) is symmetric with respect to the fixed set of this transformation: lines $\delta = \pi - 2\alpha$ and $\delta = 3\pi - 2\alpha$.

\textit{Property 6.} If $\delta = \pi - 2\alpha$ or $\delta = 3\pi - 2\alpha$, the system \eqref{ensemble} is reversible. Indeed, under these conditions, Eq. \eqref{ensemble} take the form
\begin{equation*}
\begin{cases}
\mathop{\phi_1}\limits^\cdot =\displaystyle \gamma - \sin \phi_1 + \frac{d}{1 + e^{\pm k(\sin \alpha - \sin \phi_2)}}\\
\mathop{\phi_2}\limits^\cdot =\displaystyle \gamma - \sin \phi_2 + \frac{d}{1 + e^{\pm k(\sin \alpha - \sin \phi_1)}}
\end{cases}.
\end{equation*}
Here the sign ``+'' is taken for the case $\delta = \pi - 2\alpha$. The set of points, with respect to which the phase space is symmetric, is the line $\phi_1 + \phi_2 = \pi \pmod {2\pi}$. The involution implementing this symmetry is the mapping $R: (x, y) \mapsto (\pi - y, \pi - x)$.

\textit{Property 7.} Since the coupling function %, defined by Eq.
\eqref{coupling_func} is positive, at positive values of $d$ the coordinates of equilibria in the system \eqref{ensemble} obey the inequalities $\arcsin \gamma < \phi_{1,2} < \pi - \arcsin \gamma$. 

%\textit{Property 8.}
%Phase space of the system \eqref{ensemble} does not contain limit cycles that are not symmetric to themselves. This property was proved numerically, though its origin remains unclear. We can safely state that the reason is not the presence of symmetry in the system: one can construct similar system with symmetry that possess described limit cycles.
%\begin{equation*}
%\begin{cases}
%\mathop{\phi_1}\limits^\cdot = \gamma - \sin \phi_1 + d\cdot I(\phi_2) + 0.1\arcsin \sin \phi_1\\
%\mathop{\phi_2}\limits^\cdot = \gamma - \sin \phi_2 + d\cdot I(\phi_1) + 0.1\arcsin \sin \phi_2
%\end{cases}.
%\end{equation*}

%\begin{figure}[H]
%	\centering
%	\includegraphics[width = 0.45\linewidth]{fig/PP_alpha=7pi_4-0.1_delta=3pi_4_perturbed.eps}
%	\caption{Phase portrait of perturbed system. The red lines correspond to unstable limit cycles, the green lines -- to the stable ones.}
%	\label{PP_perturbed}
%\end{figure}

\section{Dynamics of the system}
\label{sec_overall_dyn}

We have found out that the system \eqref{ensemble}, depending on the values of control parameters $\alpha$ and $\delta$ of excitatory coupling, is able to generate all main types of neuron-like activity typical for HCO: excitable steady state and in-phase/anti-phase oscillations. Below we show how these states arise and disappear in the system \eqref{ensemble} when the governing parameters are varied.

This section is organized as follows. In the first subsection we present an overall dynamical sketch of the system for the case of strong coupling. It includes, first of all, the %detailed bi
two-parameter state diagram. Then we characterize regions of multistability, proceeding to the description of the phase space and kinds of neuron-like activity for parameters taken from each diagram region. In the next paragraph the obtained states are observed in application to the HCO modelling. In the last part of the first subsection we discuss bifurcation scenarios that lead to the onset and destruction of all obtained types of neuron-like activity. The second subsection concerns evolution of the excitable state caused by variation of the coupling strength $d$. The last subsection analyzes how the variation of the coupling strength affects tonic spiking, namely, in-phase and anti-phase oscillations.

\subsection{Overall dynamical sketch for fixed coupling strength}

\begin{figure}[H]
	\centering
	(a)\includegraphics[width = 0.45\linewidth]{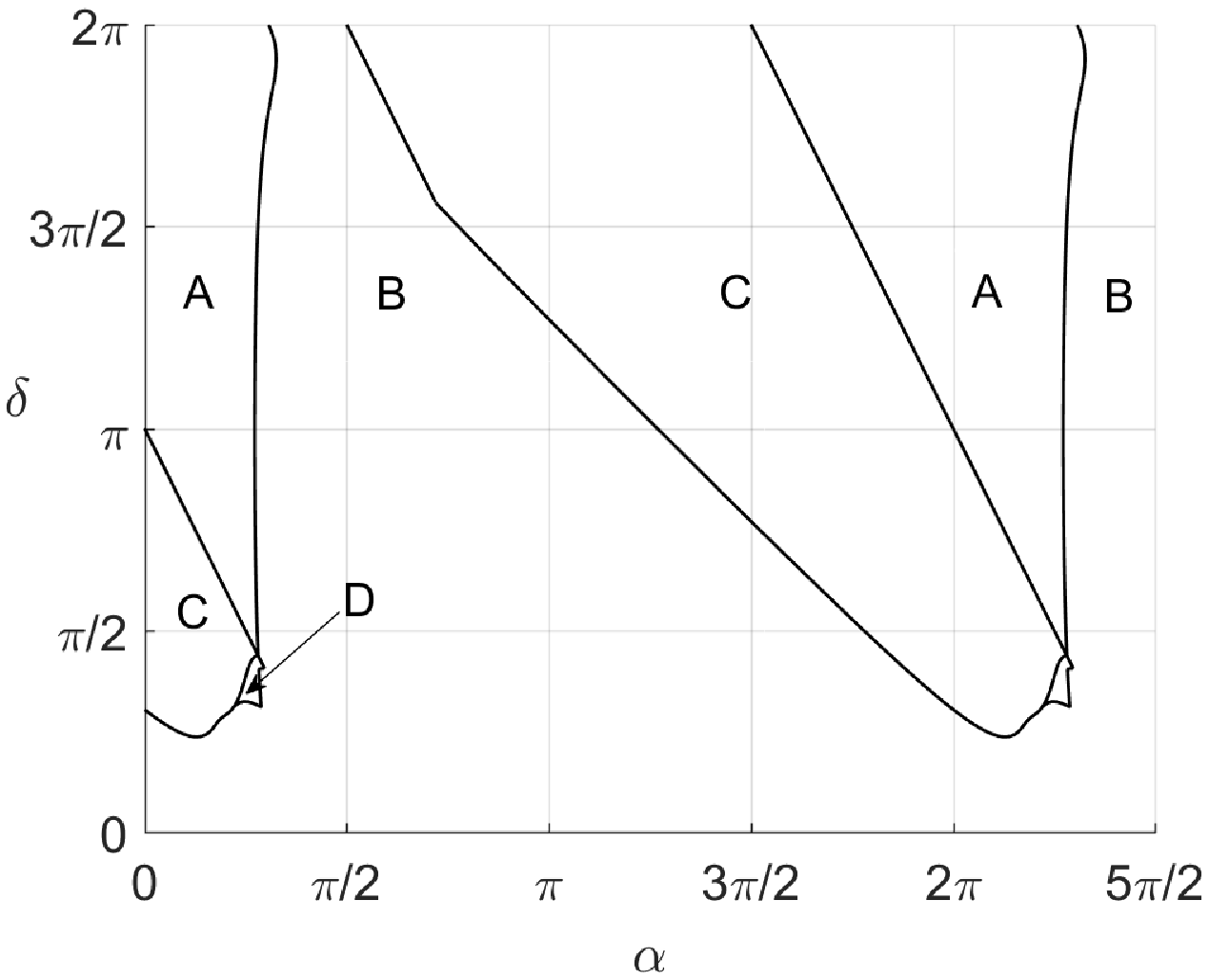}
	(b)\includegraphics[width = 0.45\linewidth]{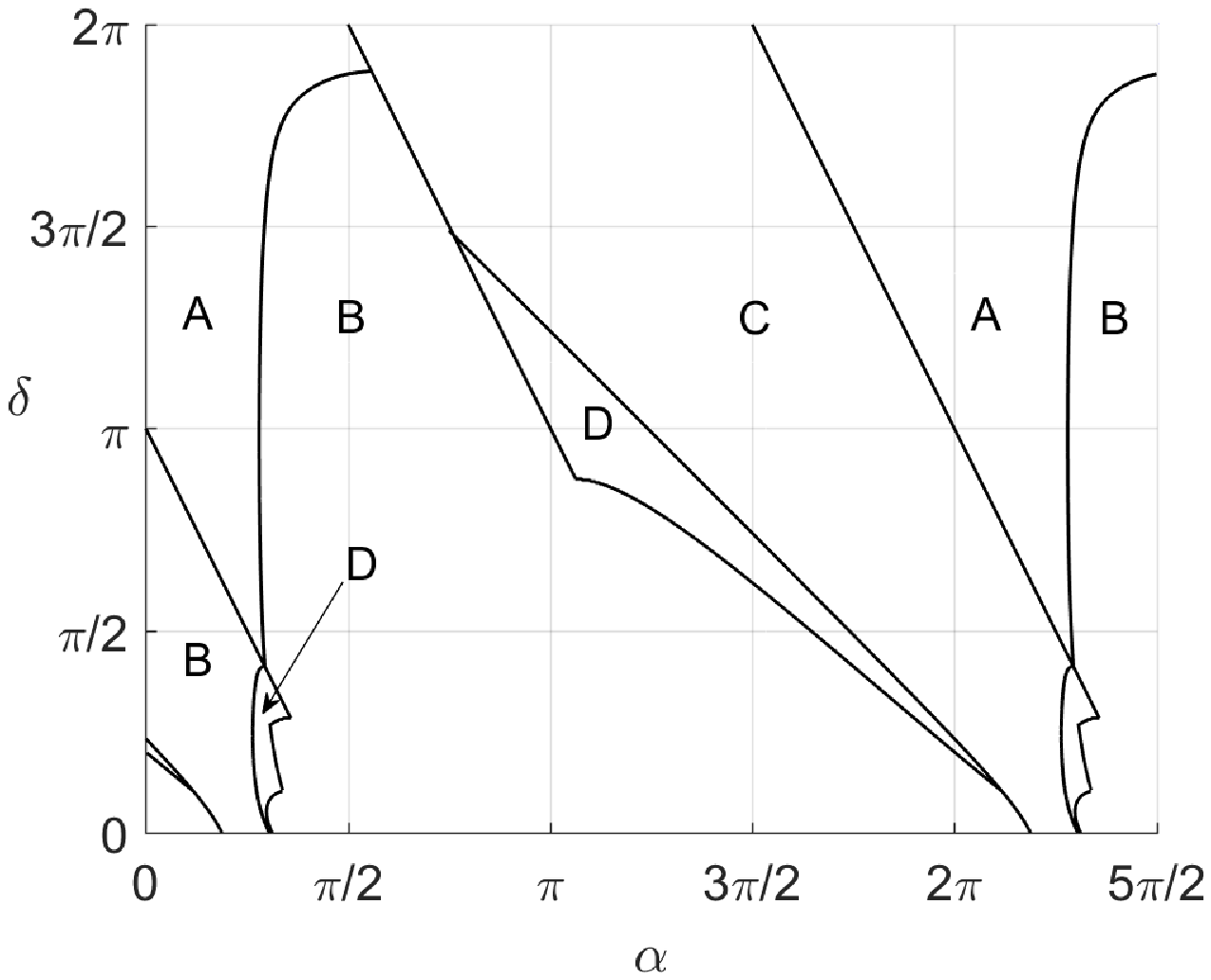}
	\caption{Map of neuron-like temporal patterns for fixed values of the coupling strength: $d=0.31$ in panel (a) and $d=1$ in panel (b). Regions correspond, respectively, to in-phase tonic spiking ($A)$, excitable state ($B$), anti-phase tonic spiking ($C$) and bistability: coexistence of the excitable state and anti-phase tonic spiking (region $D$).}
	\label{Bif_diagram}
\end{figure}

Fig.~\ref{Bif_diagram}, obtained by combining analytical and numerical methods, presents on the $(\alpha, \delta)$ parameter plane  the map of neuron-like temporal patterns. There exists a threshold value $d_{th}$ (dependent on the other system parameters) so that for the values of $d$ below $d_{th}$ the motif can exhibit only excitable behaviour, similarly to the dynamics of the single element. Increase of the coupling strength beyond $d_{th}$ leads to the onset of collective spiking dynamics. In the left panel of this Figure, the coupling strength $d$, albeit low, suffices to reproduce all main types of neuron-like behavior.  The right panel, Fig.~\ref{Bif_diagram}(b), shows locuses of different temporal patterns for the case when the value of the coupling strength $d$ is raised to $d = 1$. The main effect manifests itself in presence of the quite wide region $D$ of bistability, located between the regions $B$ (excitable state) and $C$ (anti-phase spiking). This phenomenon can be explained as follows:  in the course of increase of $d$, stability regions for the steady state and for the anti-phase limit cycle start to overlap, resulting in the coexistence of two attractors in the phase space. The borderlines of other regions of neuron-like temporal patterns are also shifted when $d$ is increased, and at appropriate values of $\alpha$ and $\delta$ the excitable state gets replaced by oscillatory activity (both in-phase and anti-phase).

Let us list the types of neuron-like activity, observable in each of the regions from Fig. \ref{Bif_diagram}.

The region $A$ features in-phase spiking activity with $\phi_1(t)=\phi_2(t)$. In the phase space, mathematical image of this activity type is the stable in-phase limit cycle. In the region $B$ only the excitable state exists. Although dynamics in $B$ is simple, it corresponds to different stable equilibria. From the point of view of neuroscience, coexistence of different excitable states can describe different conditions of the membrane potential of neuron-like elements, including depolarization and hyperpolarization. In the region $C$ the system \eqref{ensemble} exhibits only anti-phase spiking activity, described by the stable anti-phase limit cycle. The region $D$ is the only domain of bistability, where anti-phase spiking patterns coexist with excitable behavior. 

\begin{figure}[H]
	\centering
	(a)\includegraphics[width = 0.21\linewidth]{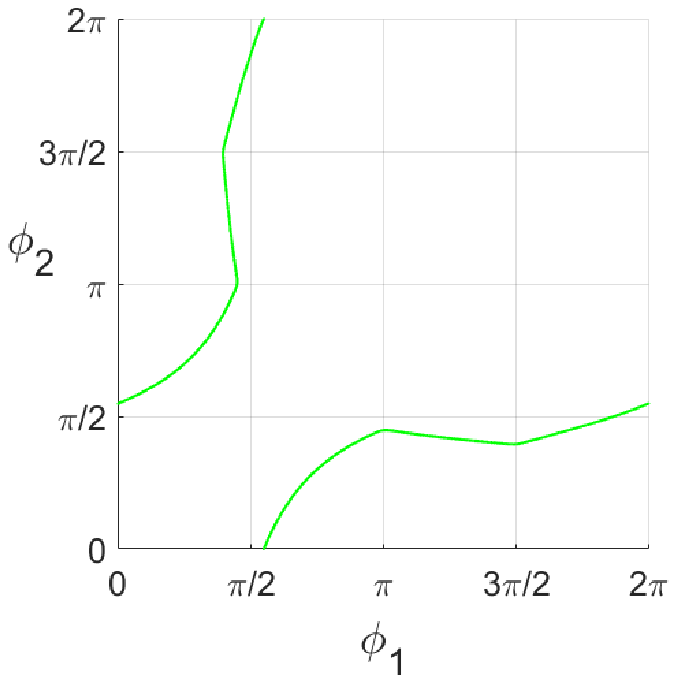}
	(b)\includegraphics[width = 0.21\linewidth]{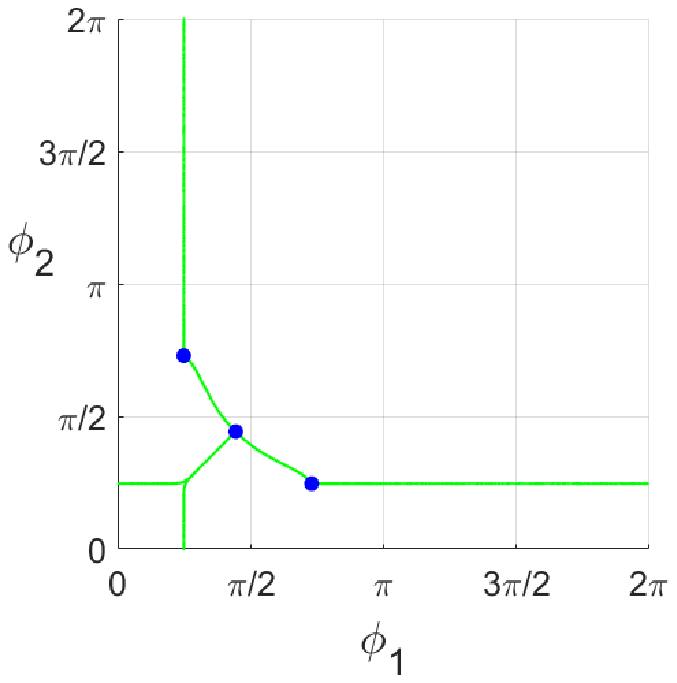}
	(c)\includegraphics[width = 0.21\linewidth]{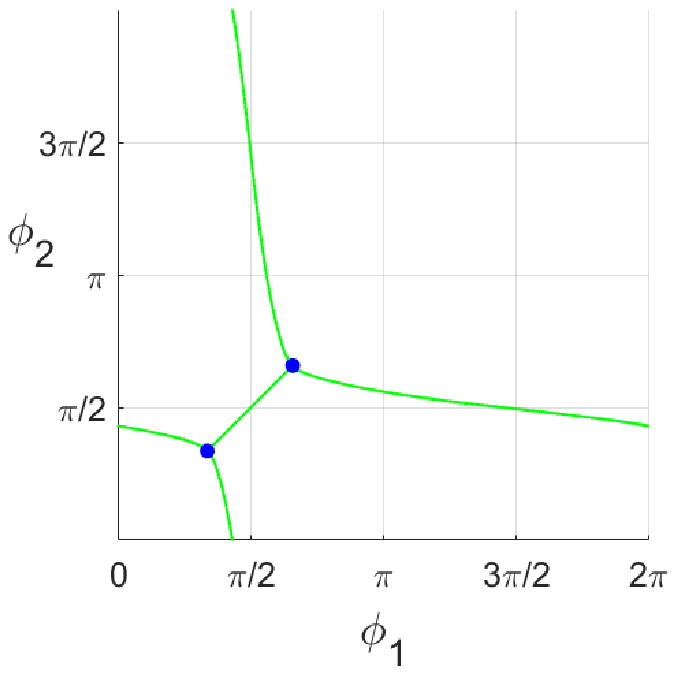}
	(d)\includegraphics[width = 0.21\linewidth]{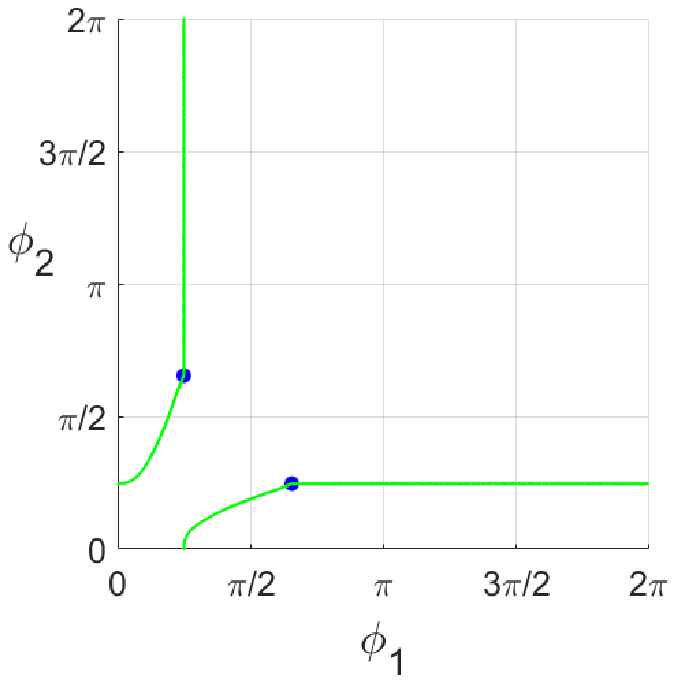}\\
	(e)\includegraphics[width = 0.2\linewidth]{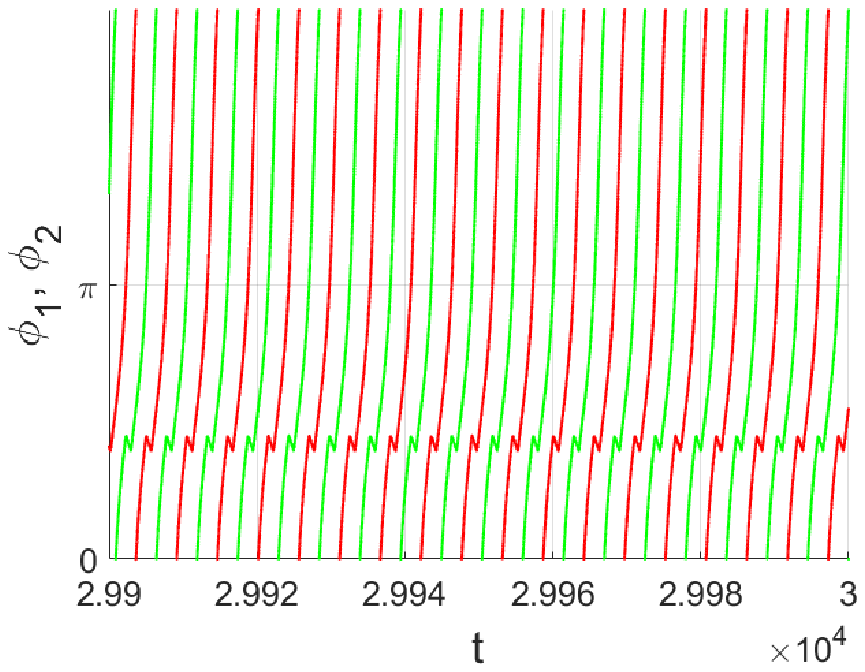}
	(f)\includegraphics[width = 0.2\linewidth]{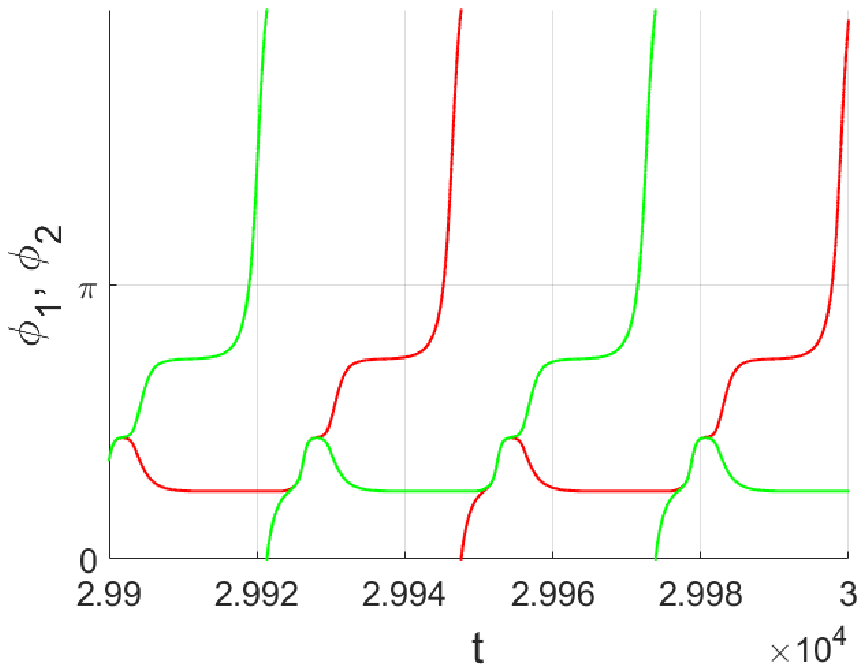}
	(g)\includegraphics[width = 0.2\linewidth]{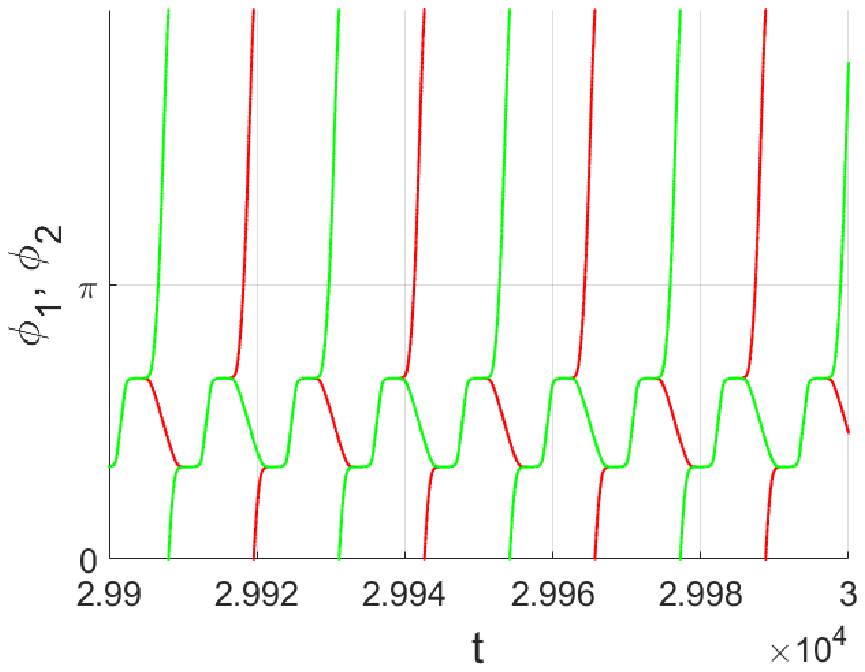}
	(h)\includegraphics[width = 0.2\linewidth]{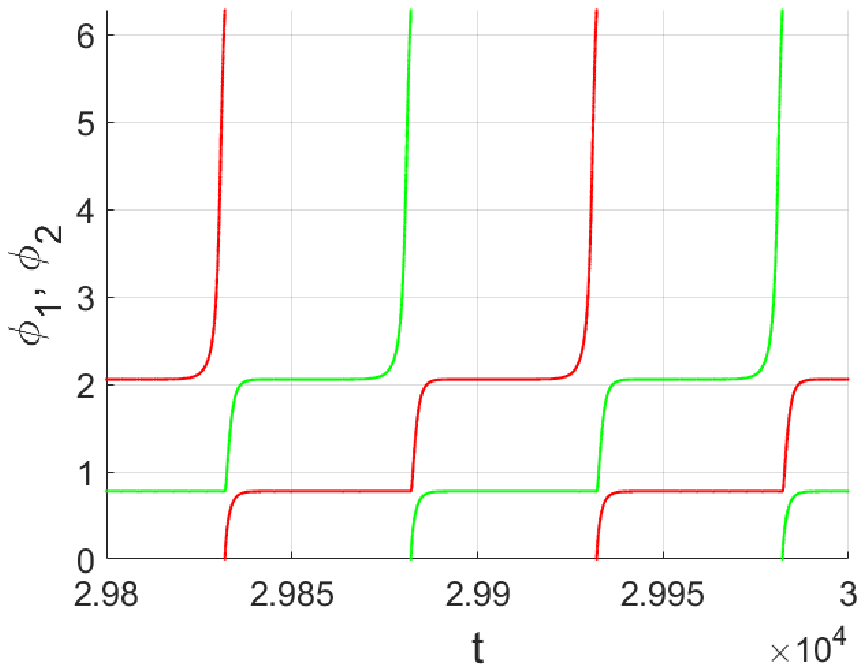}
	\caption{Examples of anti-phase spiking states: phase space (top row) and time series (bottom). Blue dots denote saddles. (a) $\alpha = 3\pi/2$, $\delta = 3\pi/2$. (b) $\alpha = 1$, $\delta = 0.296$. (c) $\alpha = 1.124$, $\delta = 0.8755298$. (d) $\alpha = 0.2$, $\delta = 0.4707920318$. (e) $\alpha = 3\pi/2$, $\delta = 3\pi/2$. (f) $\alpha = 1$, $\delta = 0.296$. (g) $\alpha = 1.124$, $\delta = 0.8755298$. (h) $\alpha = 0.2$, $\delta = 0.4707920318$.}
	\label{Time_series}
\end{figure}

In the framework of HCO modelling the most interesting and %valuable 
important states are those with anti-phase activity. Fig.~\ref{Time_series} renders time series of stable anti-phase limit cycles along with their images in the phase space at different values of governing parameters.

Let us have a closer look at transitions between observed types of temporal patterns of neuron-like activity in Eqs. \eqref{ensemble}. We restrict ourselves to $d=1$.

\begin{figure}[H]
  \centering
  (a)\includegraphics[width = 0.3\linewidth]{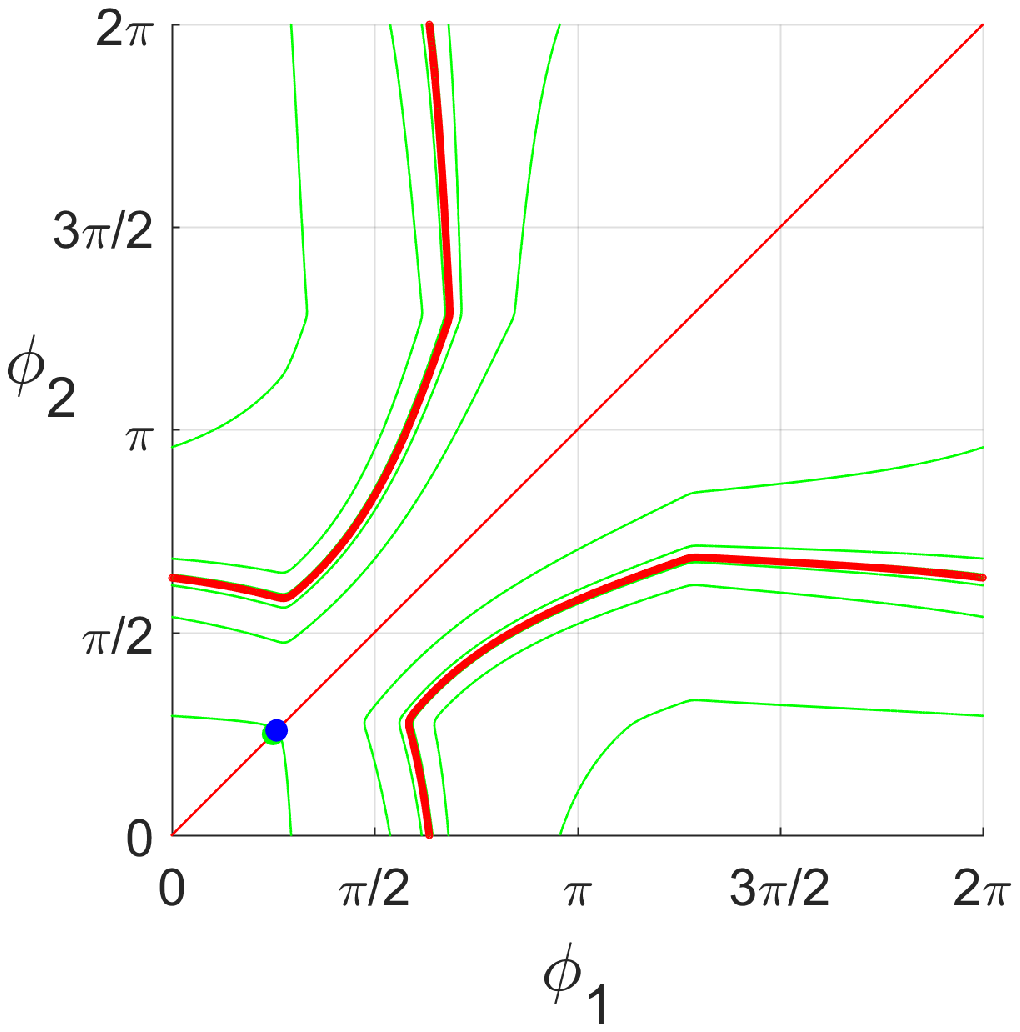}
  (b)\includegraphics[width = 0.3\linewidth]{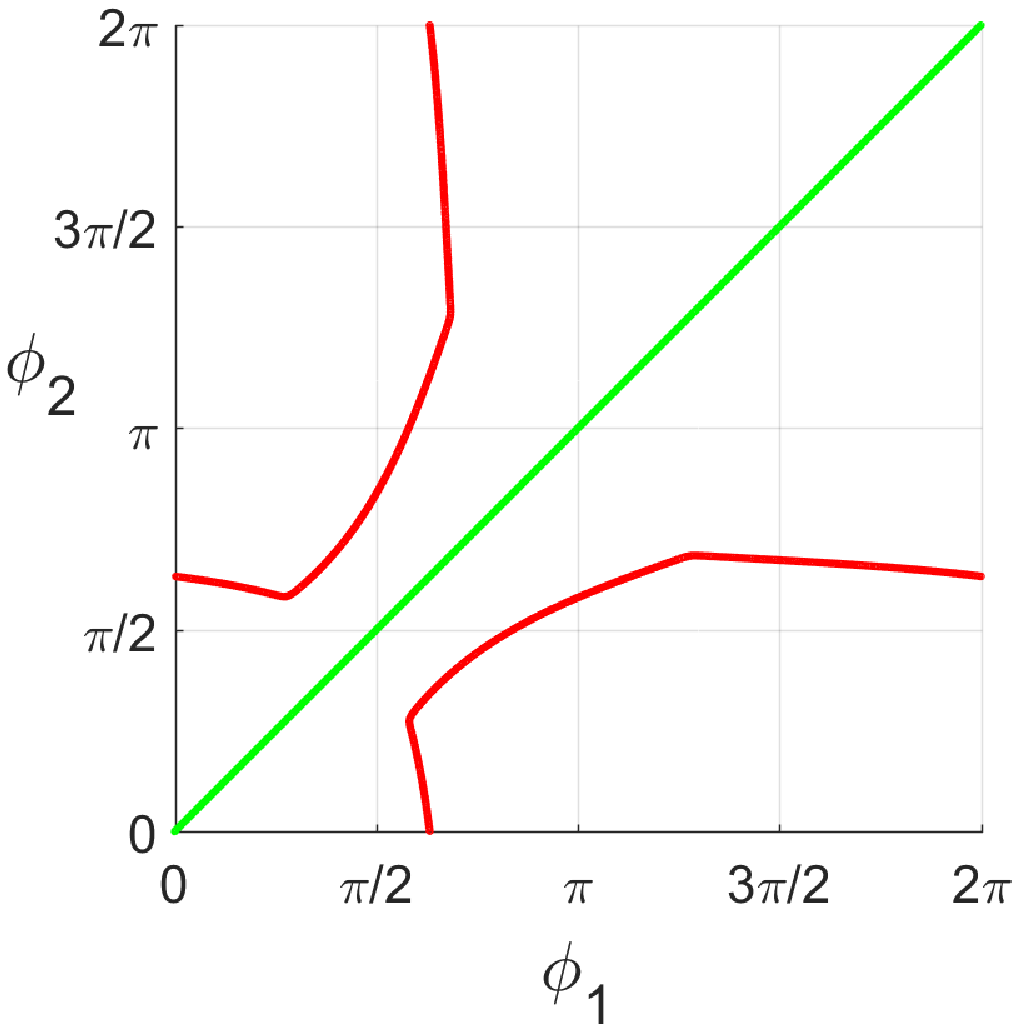}
  \caption{Scenario of birth of the in-phase limit cycle at the borderline between the regions $B$ and $A$. $\delta = \pi$. (a) $\alpha = 0.885$, (b) $\alpha = 0.875$. Here red/green bold lines correspond to the unstable/stable limit cycle, red/green curves correspond to unstable/stable separatrices. Green dot corresponds to stable equilibrium (stable node), blue dot -- to saddle equilibrium.}
  \label{Born_inphase}
\end{figure}

We start with the transition between the regions $B$ and $A$. To this end, we fix $\delta = \pi$ and decrease the governing parameter $\alpha$ from the value $\alpha = 0.885$ to the value $\alpha = 0.875$, crossing thereby the borderline between these regions (see Fig. \ref{Born_inphase}). As a result of the saddle-node bifurcation on the invariant curve, the stable in-phase limit cycle appears in the phase space.

\begin{figure}[H]
  \centering
  (a)\includegraphics[width = 0.3\linewidth]{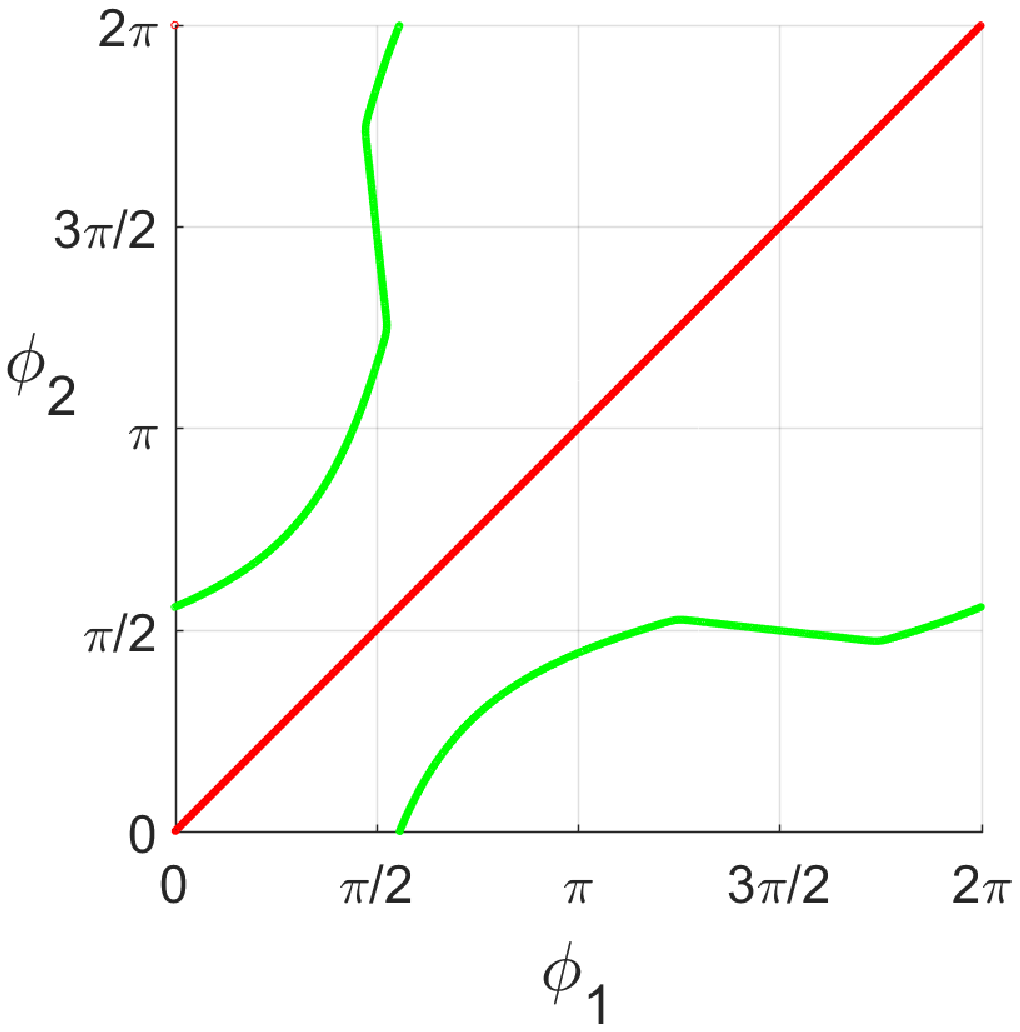}
  (b)\includegraphics[width = 0.3\linewidth]{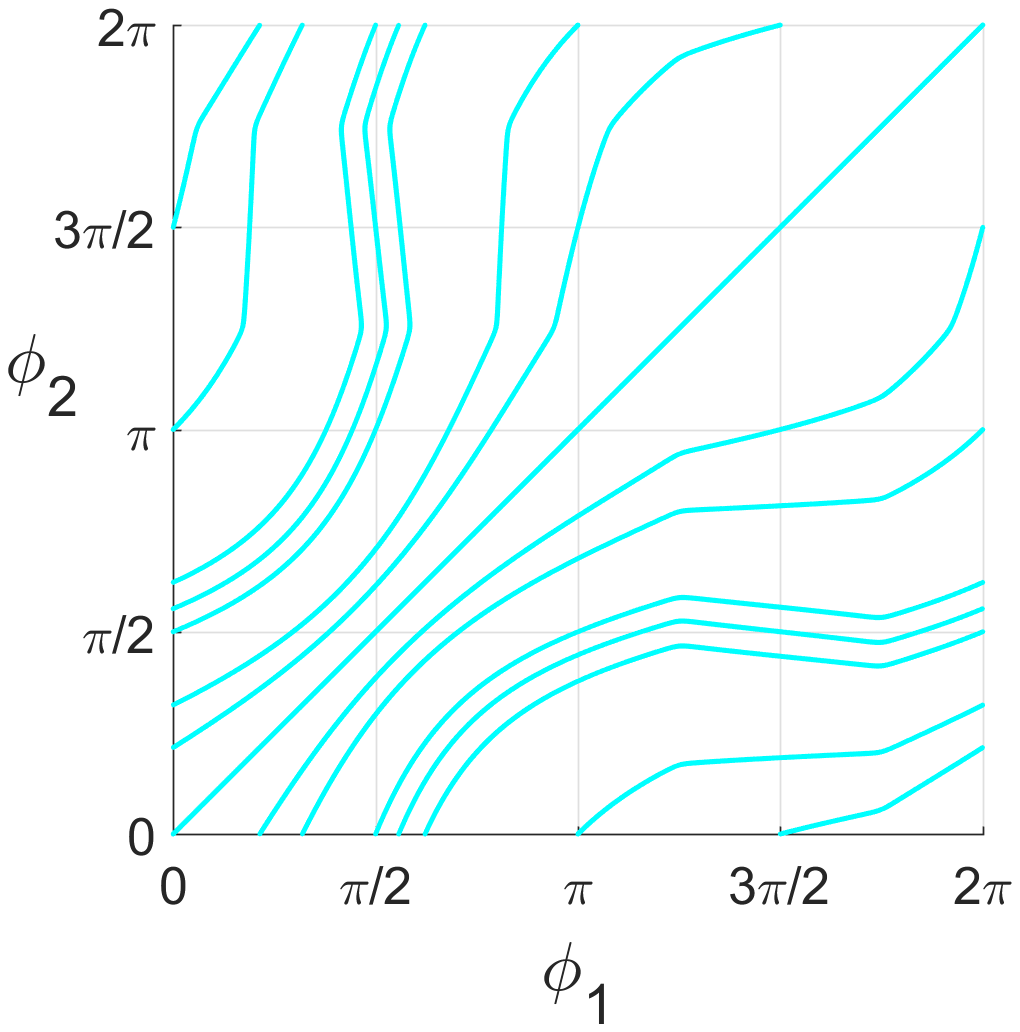}
  (c)\includegraphics[width = 0.3\linewidth]{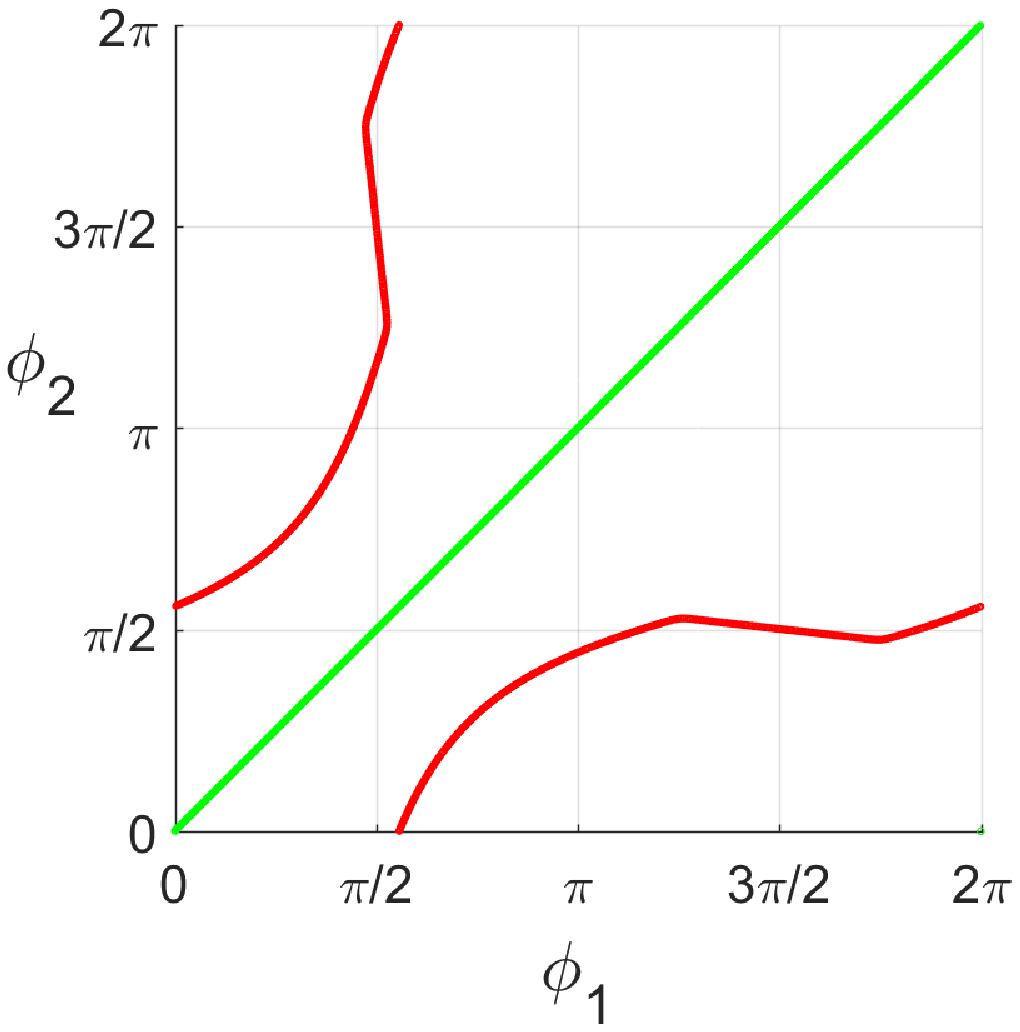}
  \caption{Birth scenario for the in-phase limit cycle at the borderline between regions $C$ and $A$. $\delta =\displaystyle\frac{3\pi}{2}$. (a) $\alpha = \displaystyle\frac{7\pi}{4} - 0.01$. (b) $\alpha = \displaystyle\frac{7\pi}{4}$. (c) $\alpha = \displaystyle\frac{7\pi}{4} + 0.01$. In (a) red/green curves correspond to the unstable/stable in-phase/anti-phase cycle. In (b) light blue curves mark closed trajectories that pass through each point of the phase space. In (c) green/red curves correspond to the stable/unstable in-phase/anti-phase cycles.}
  \label{Born_inphase1}
\end{figure}

The transition from region $C$ to the region $A$ is more involved. To illustrate the pertinent bifurcation scenario we fix $\delta = \frac{3\pi}{2}$ and build phase portraits of the system for values of parameter $\alpha$ taken from the region $C$ before the transition, on the borderline between two regions and in the region $A$ right after the bifurcation. Fig.~\ref{Born_inphase1} shows the bifurcation, as a result of which the in-phase limit cycle becomes stable. In Fig.~\ref{Born_inphase1}(a), an unstable in-phase and a stable anti-phase cycles are present. When a parameter $\alpha$ reaches its bifurcation value $\alpha=\frac{7\pi}{4}$ (see Fig. \ref{Born_inphase1}(b)), a continuum of closed trajectories exists; %these orbits pass through each point of the phase space. 
the whole torus is foliated into neutrally stable periodic orbits.
In the course of this non-local bifurcation, the in-phase cycle acquires stability, whereas the anti-phase one gets destabilized (Fig.~\ref{Born_inphase1}(c)).

\begin{figure}[H]
  \centering
  (a)\includegraphics[width = 0.3\linewidth]{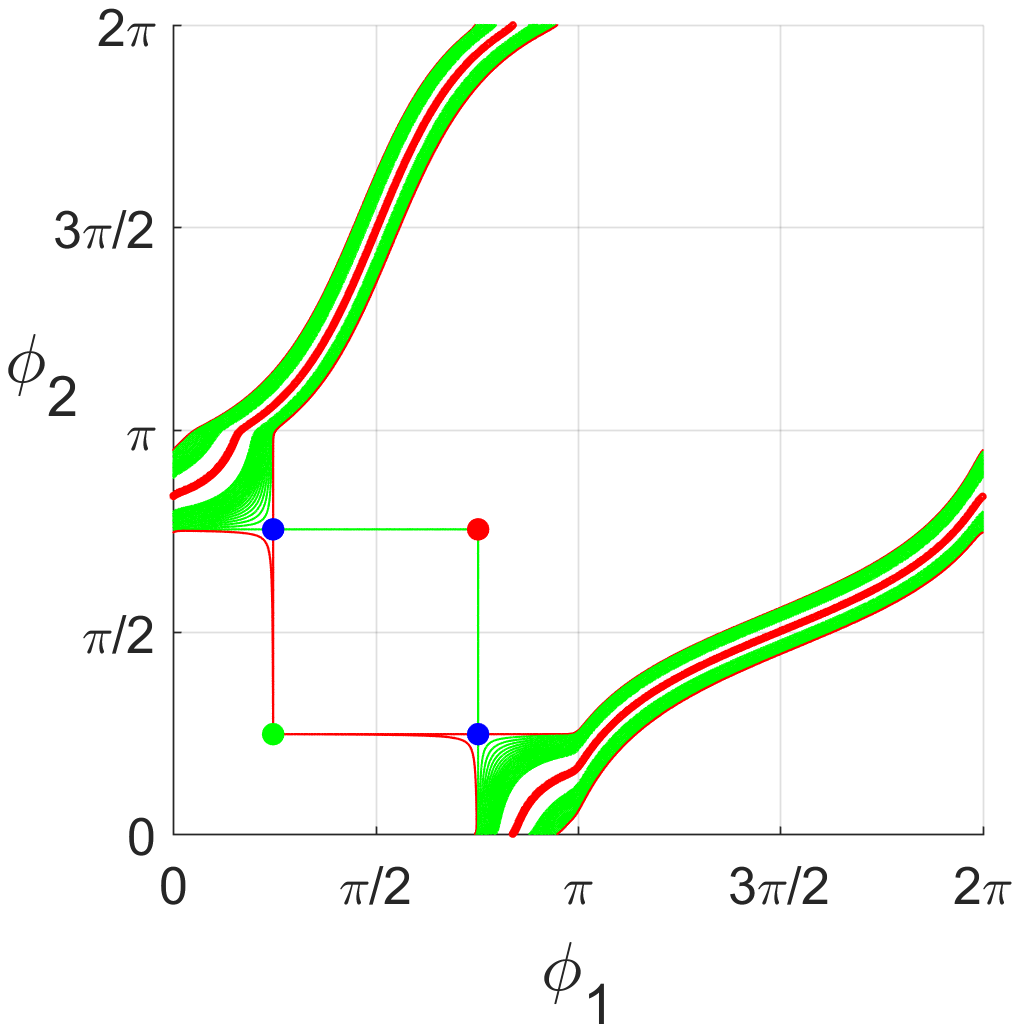}
  (b)\includegraphics[width = 0.3\linewidth]{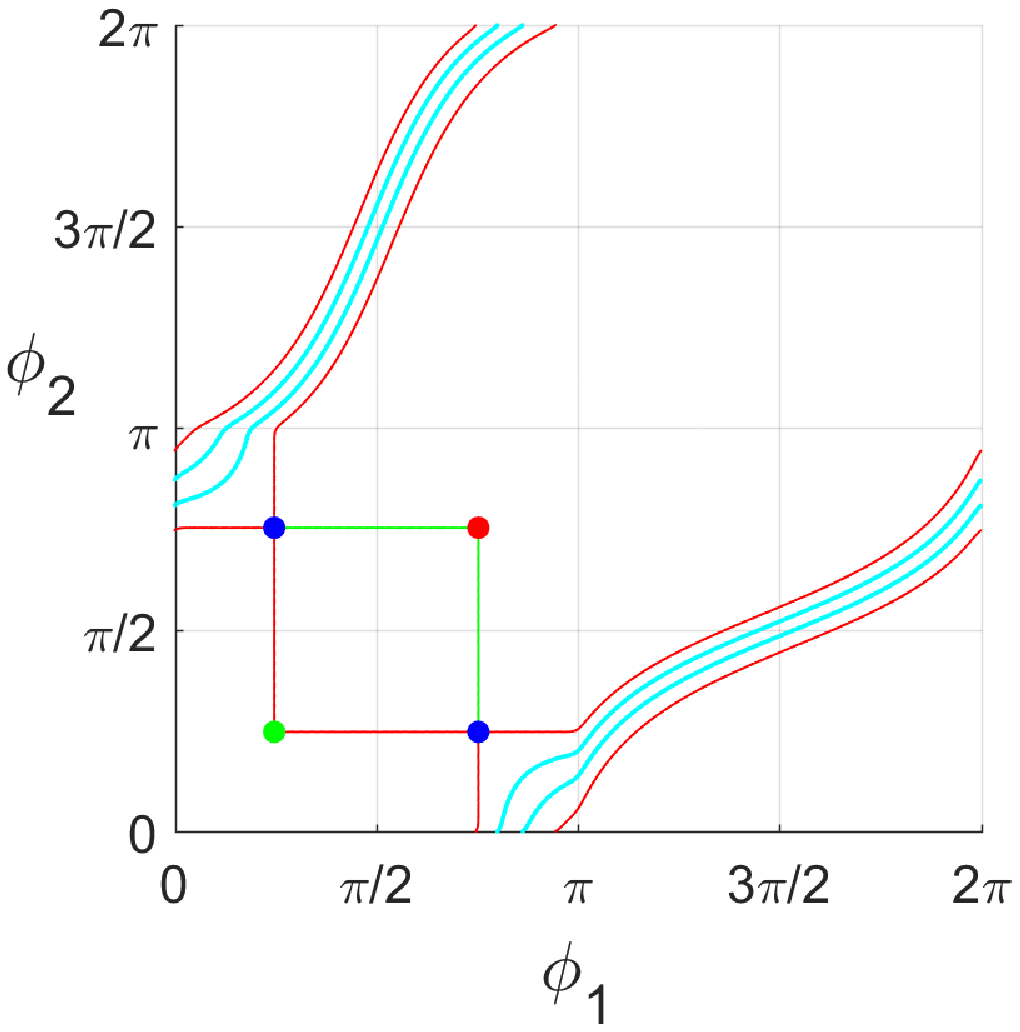}
  (c)\includegraphics[width = 0.3\linewidth]{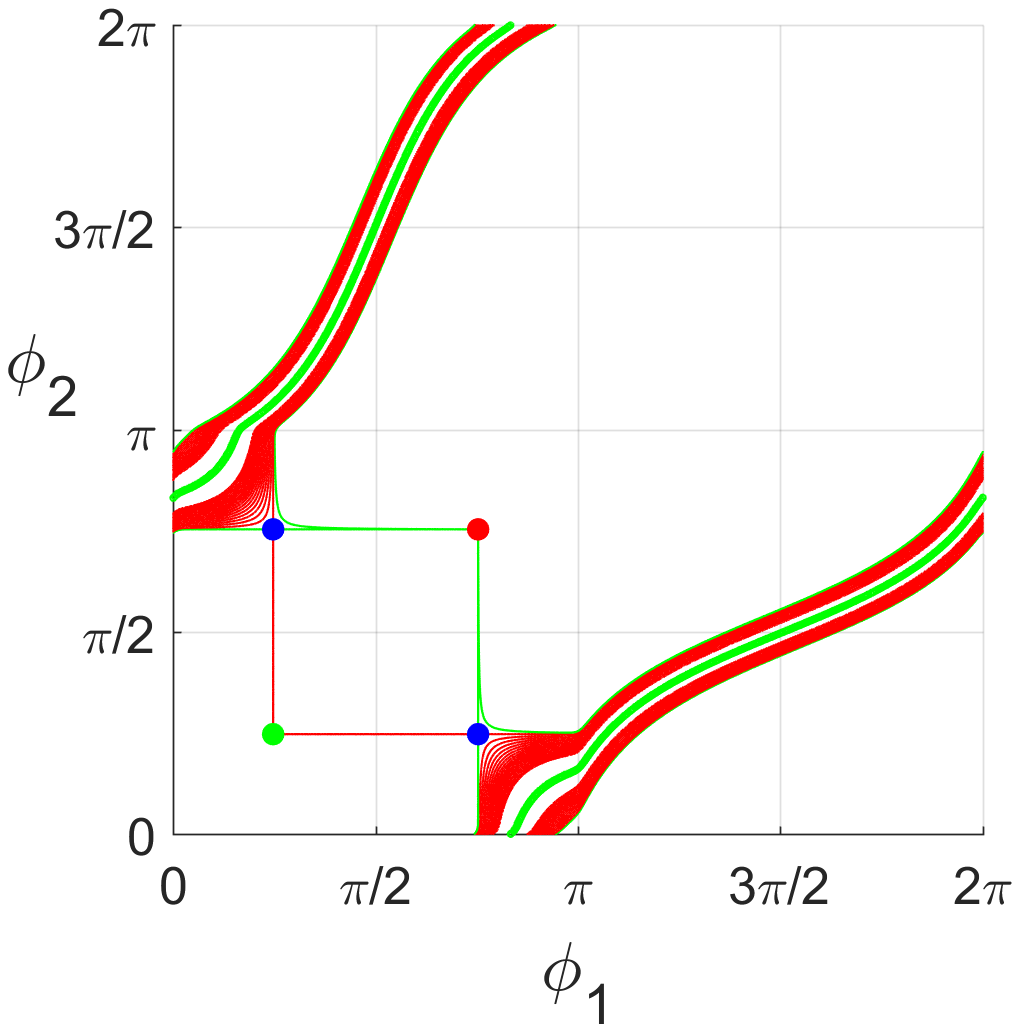}
  \caption{The first birth scenario of the anti-phase limit cycle at the borderline between the regions $C$ and $D$. $\delta = \pi$. (a) $\alpha = \pi - 0.01$. (b) $\alpha = \pi$. (c) $\alpha = \pi + 0.01$. Here red/green curves correspond to unstable/stable separatrices, red dot corresponds to unstable equilibrium (unstable node), green dot -- to stable equilibrium (stable node), blue dot -- to saddle equilibrium. In (a) red bold line corresponds to the unstable anti-phase cycle. In (b) light blue curves mark closed trajectories that pass through each point of the area of the phase space, bounded by homoclinic trajectories. In (c) the green bold line corresponds to the stable anti-phase cycle.}
  \label{Born_antiphase1}
\end{figure}

The sophisticated borderline between the regions $C$ and $D$ offers several scenarios of the birth of bistability between anti-phase spiking pattern and the excitable state. The first scenario is presented in Fig. \ref{Born_antiphase1}. For $\alpha = \pi - 0.01$, the unstable anti-phase limit cycle exists in the phase space, so that unstable saddle separatrices tend to the stable state of equilibrium. One stable separatrix of each saddle begins at the unstable equilibrium, and the other two come from the unstable limit cycle. In the course of the bifurcation ($\alpha = \pi$), two homoclinic trajectories are formed: they delineate the region of the phase space, inside which the closed trajectories pass through each point. At $\alpha = \pi + 0.01$ the stable anti-phase limit cycle exists in the phase space. Stable saddle separatrices now begin at the unstable equilibrium. One of unstable separatrices of each saddle leads to the stable equilibrium, the other two are attracted by the stable limit cycle. 

\begin{figure}[H]
  \centering
  (a)\includegraphics[width = 0.3\linewidth]{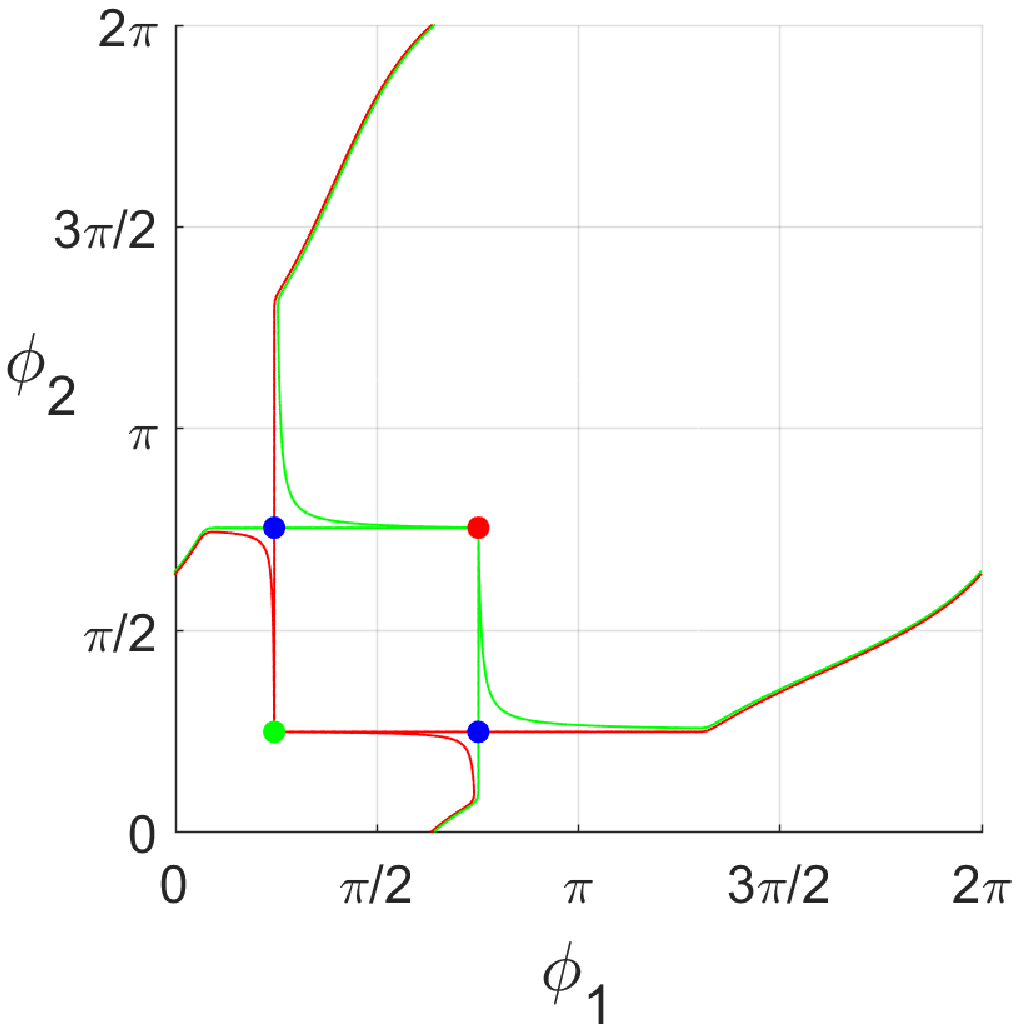}
  (b)\includegraphics[width = 0.3\linewidth]{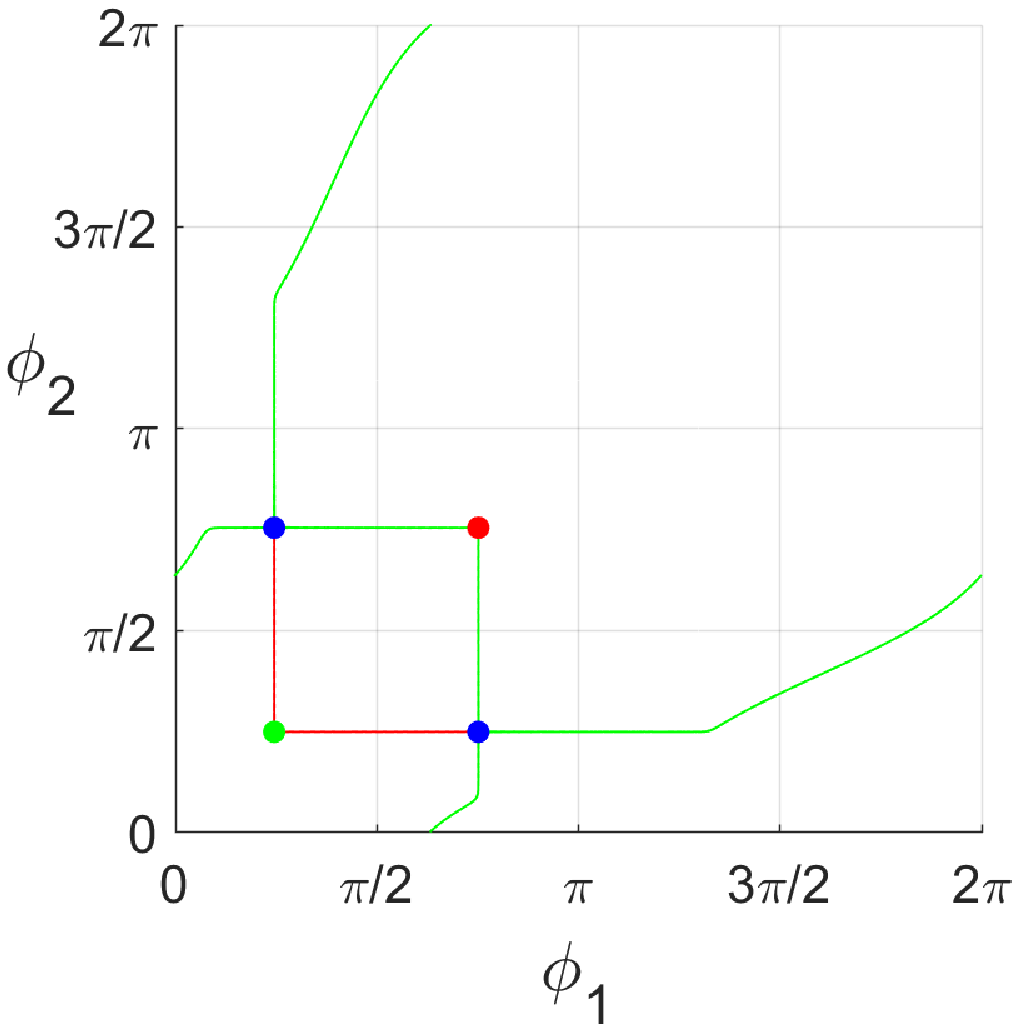}
  (c)\includegraphics[width = 0.3\linewidth]{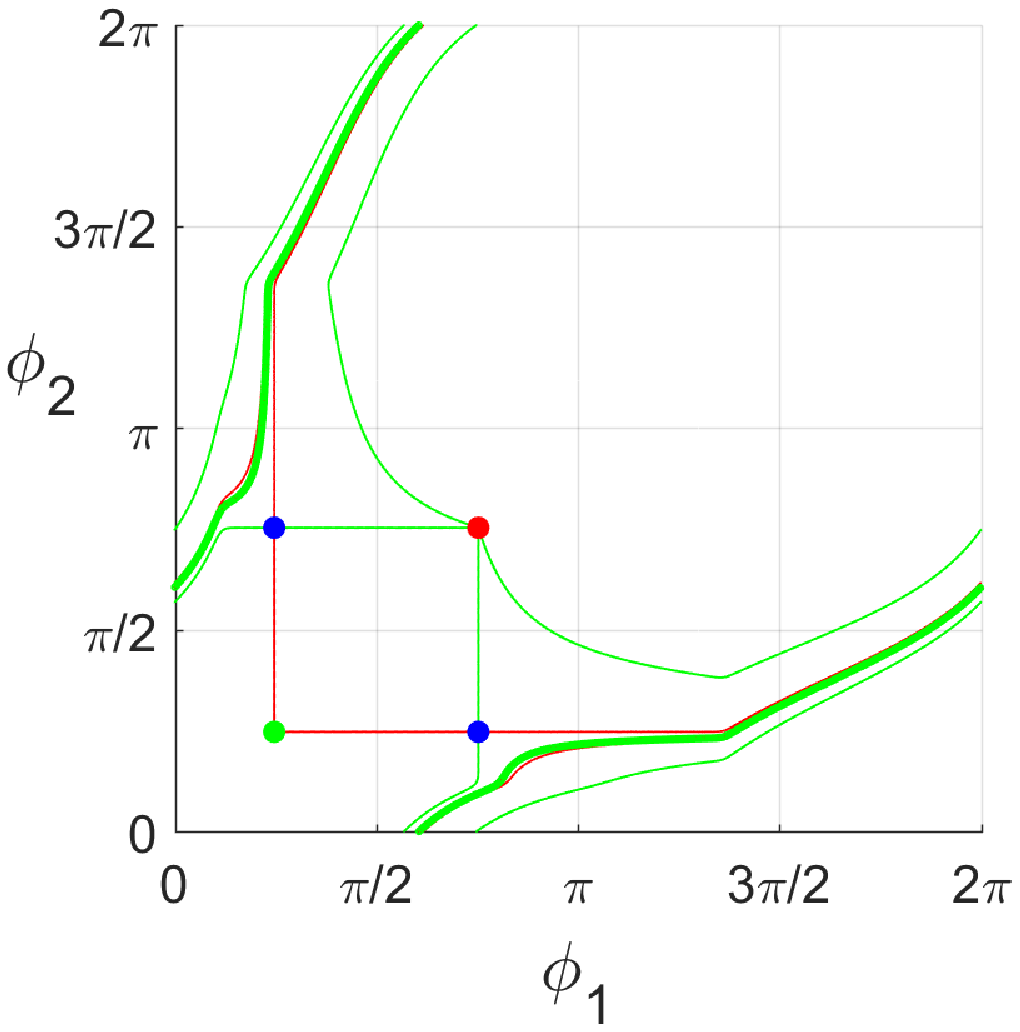}
  \caption{The second birth scenario of the anti-phase limit cycle at the borderline between the regions $C$ and $D$. $\delta =\displaystyle\frac{3\pi}{4}$. (a) $\alpha = 4.15$. (b) $\alpha = 4.1691$. (c) $\alpha = 4.28$. Here red/green curves correspond to unstable/stable separatrices, red dot corresponds to unstable equilibrium (unstable node), green dot -- to stable equilibrium (stable node), blue dot -- to saddle equilibrium. In (c) the green bold line corresponds to the stable anti-phase cycle.}
  \label{Born_antiphase2}
\end{figure}

The second scenario of birth of the stable anti-phase limit cycle in the course of transition from the region $C$ to the region $D$ is illustrated in Fig.~\ref{Born_antiphase2} and involves formation of the heteroclinic cycle (Fig. \ref{Born_antiphase2}(b)). In Figure \ref{Born_antiphase2}(a) all unstable separatrices of the saddle tend to the stable equilibrium. If we continue to increase the value of $\alpha$ up to $\alpha_{bif} \approx 4.1691$, a pair of heteroclinic trajectories between two saddles is formed in the phase space. These heteroclinic trajectories, together with the saddles, comprise a heteroclinic cycle shown in Fig.~\ref{Born_antiphase2}(b). After the bifurcation, the stable anti-phase limit cycle which attracts two unstable separatrices of the saddles, branches off the heteroclinic cycle, see Fig.~\ref{Born_antiphase2}(c).

\begin{figure}[H]
  \centering
  (a)\includegraphics[width = 0.3\linewidth]{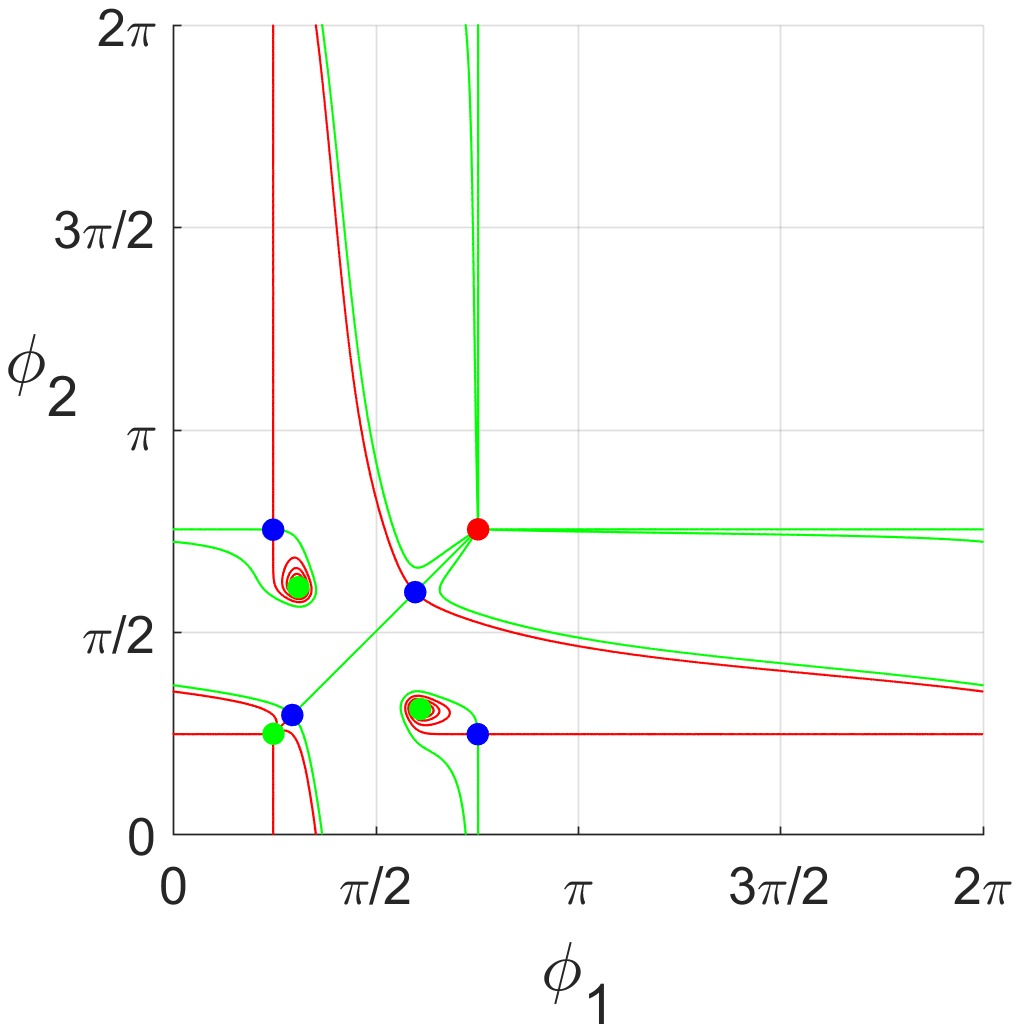}
  (b)\includegraphics[width = 0.3\linewidth]{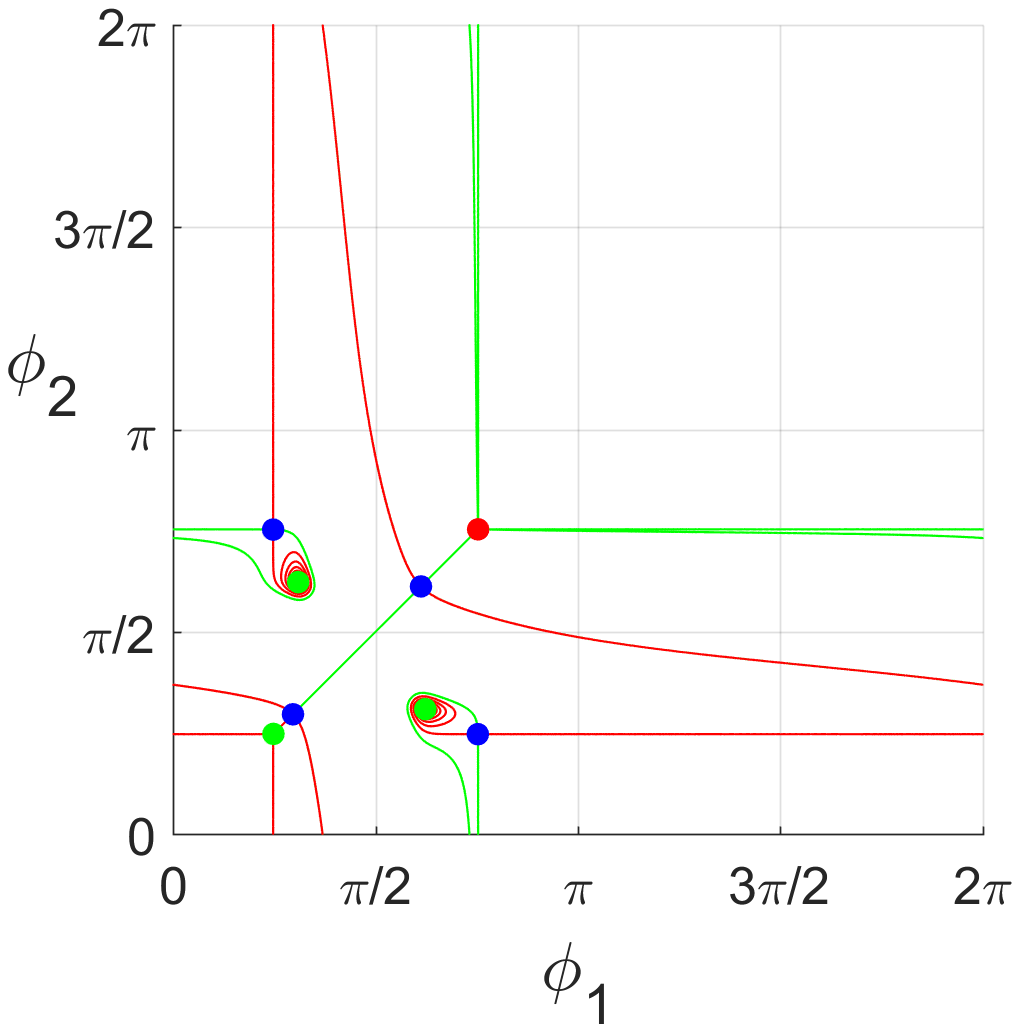}
  (c)\includegraphics[width = 0.3\linewidth]{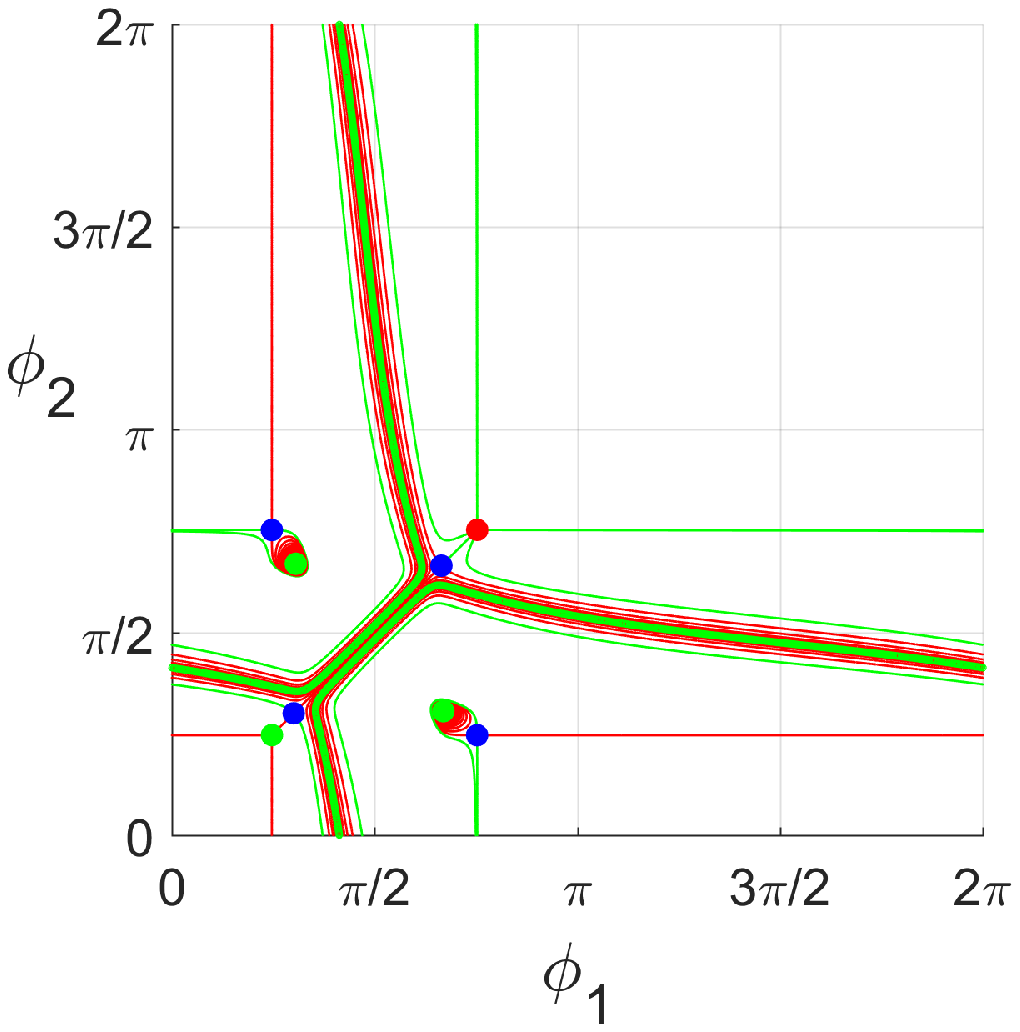}
  \caption{The third birth scenario of anti-phase limit cycle at the borderline between regions $C$ and $D$. $\alpha = 1.026$. (a) $\delta = 0.8$. (b) $\delta = 0.8432$. (c) $\delta = 1$. Here red/green curves correspond to unstable/stable separatrices, red dot corresponds to unstable equilibrium (unstable node), green dot -- to stable equilibrium (stable node), blue dot -- to saddle equilibrium. In (c) the green bold line corresponds to the stable anti-phase cycle.}
  \label{Born_antiphase3}
\end{figure}

The third scenario is observed e.g. if $\alpha$ is fixed at $1.026$ and the value of $\delta$ is raised from $\delta = 0.8$ up to $\delta = 1$. In Fig. \ref{Born_antiphase3}(a), there exists one heteroclinic trajectory lying on a diagonal line between two diagonal saddles. All other unstable separatrices as $t \rightarrow \infty$ tend to stable equilibria. Fig.\ref{Born_antiphase3}(b) shows two more symmetric heteroclinic trajectories between the same saddles. These two heteroclinic trajectories, as well as two diagonal saddles and a heteroclinic trajectory between them, lying on the diagonal, form a heteroclinic cycle. Figure Fig.\ref{Born_antiphase3}(c) shows the stable anti-phase limit cycle that appeared from heteroclinic cycle.

\begin{figure}[H]
  \centering
  (a) \includegraphics[width = 0.3\linewidth]{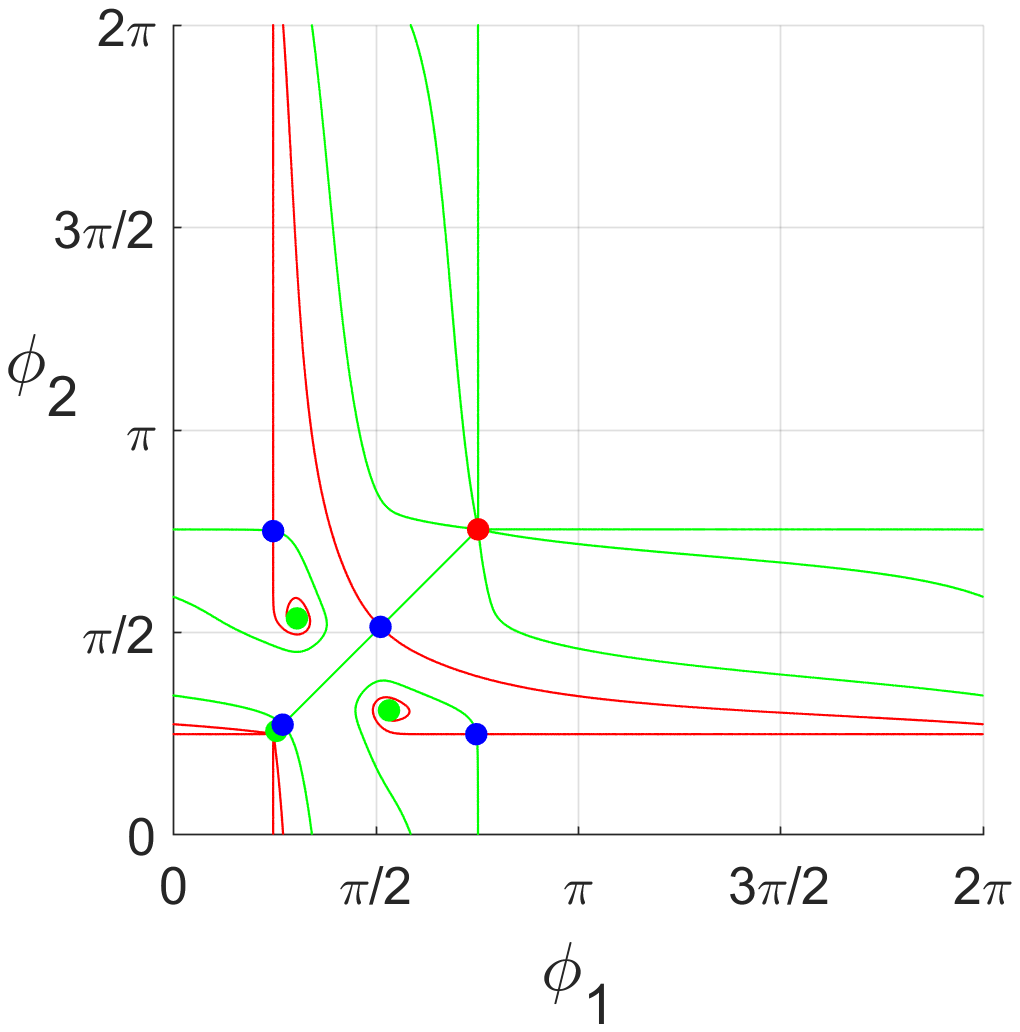}
  (b) \includegraphics[width = 0.3\linewidth]{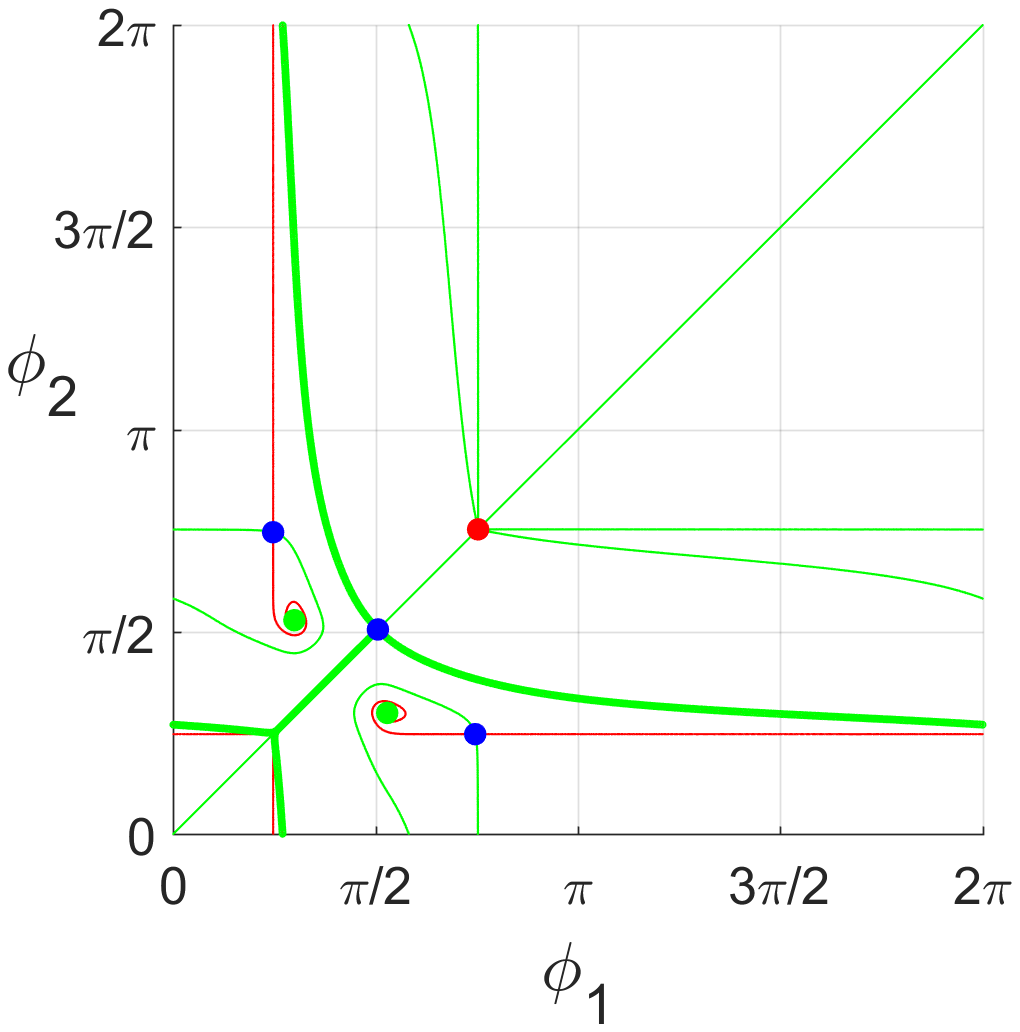}
  \caption{The fourth birth scenario of anti-phase limit cycle at the borderline between regions $C$ and $D$. $\delta = 0.53$. (a) $\alpha = 1.02$. (b) $\alpha = 1$. Here red/green curves correspond to unstable/stable separatrices, red dot corresponds to unstable equilibrium (unstable node), green dot -- to stable equilibrium (stable node), blue dot -- to saddle equilibrium. In (b) the green bold line corresponds to the stable anti-phase cycle.}
  \label{Born_antiphase4}
\end{figure}

The fourth scenario evolves as follows. Right before the bifurcation, an invariant curve exists in the phase space. It contains two saddle points, one stable equilibrium on the diagonal and the separatrices that connect them (Fig.~\ref{Born_antiphase4}(a)). On this invariant curve the saddle-node bifurcation takes place, and, as a result, the stable anti-phase limit cycle emerges (Fig.~\ref{Born_antiphase4}(b)). 

\begin{figure}[H]
  \centering
  (a) \includegraphics[width = 0.3\linewidth]{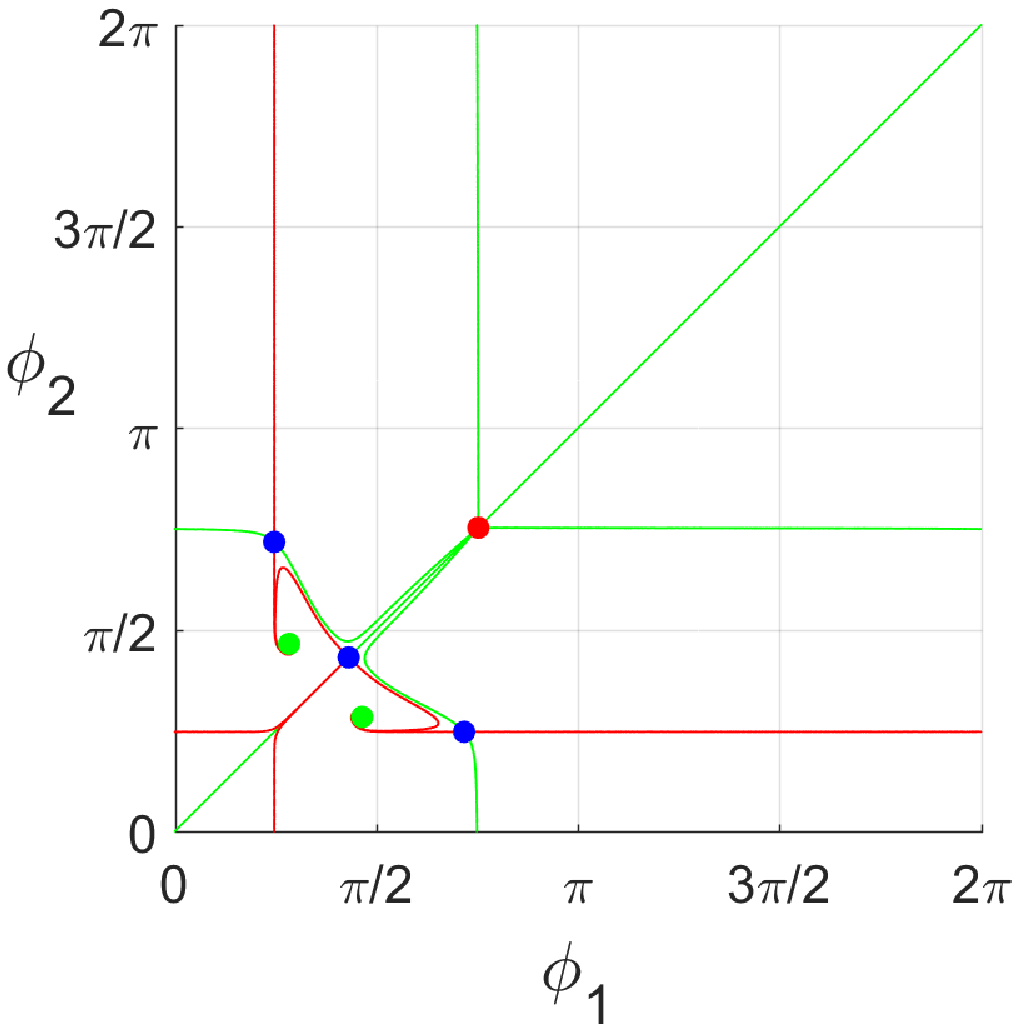}
  (b) \includegraphics[width = 0.3\linewidth]{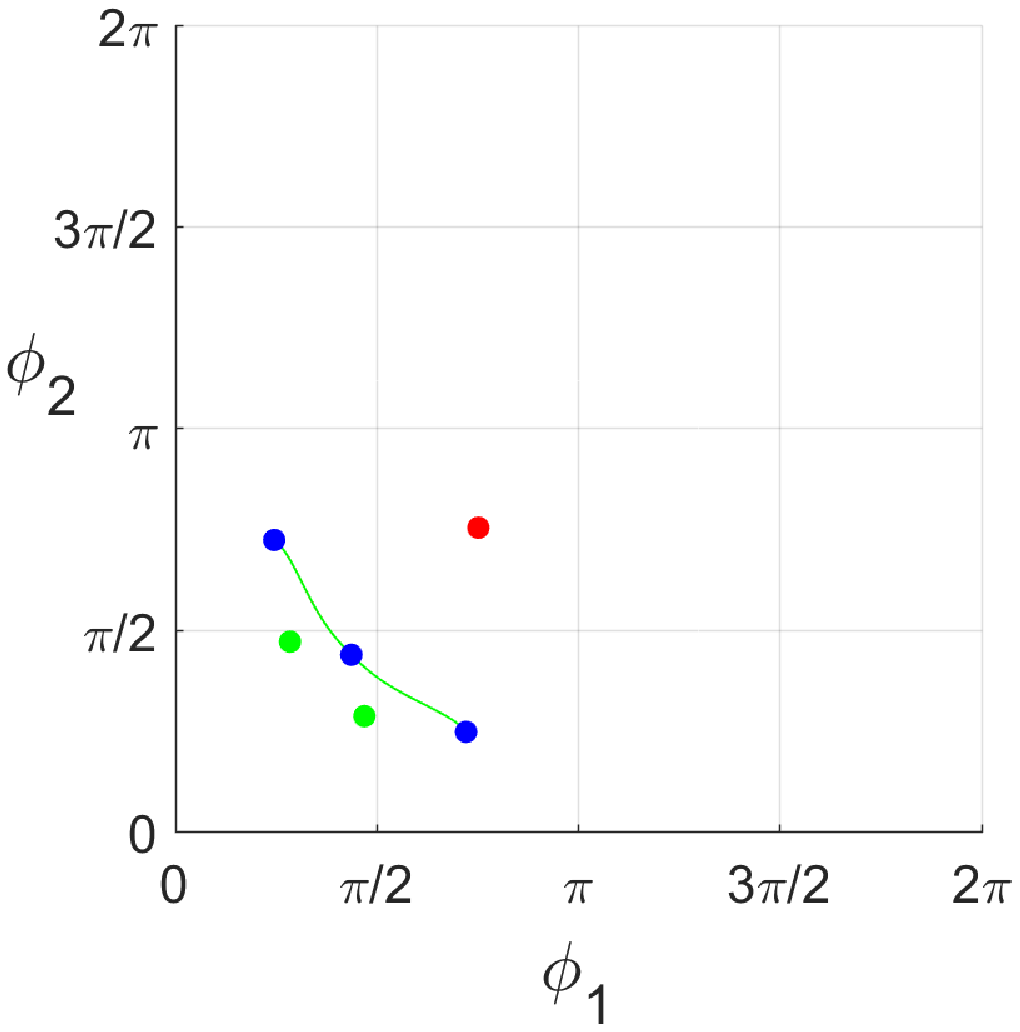}\\
  (c) \includegraphics[width = 0.3\linewidth]{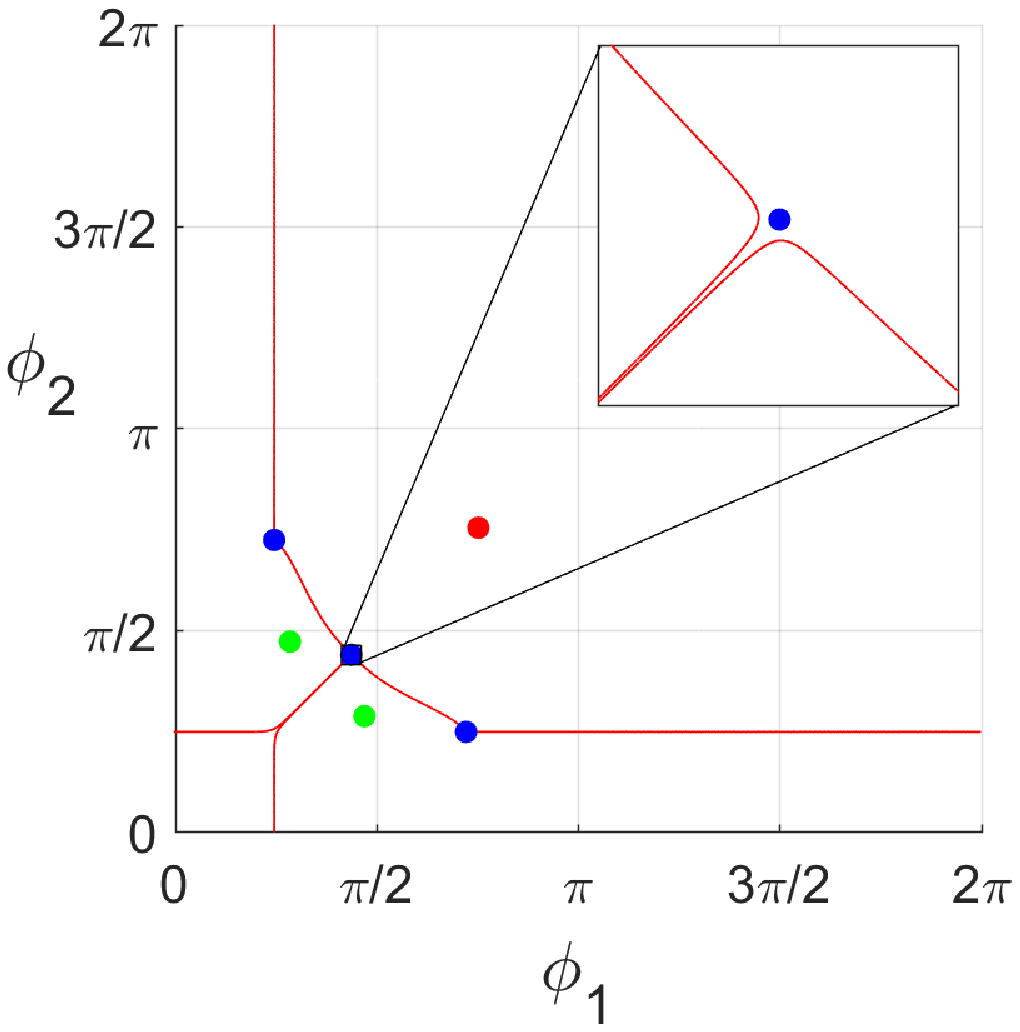}
  (d) \includegraphics[width = 0.3\linewidth]{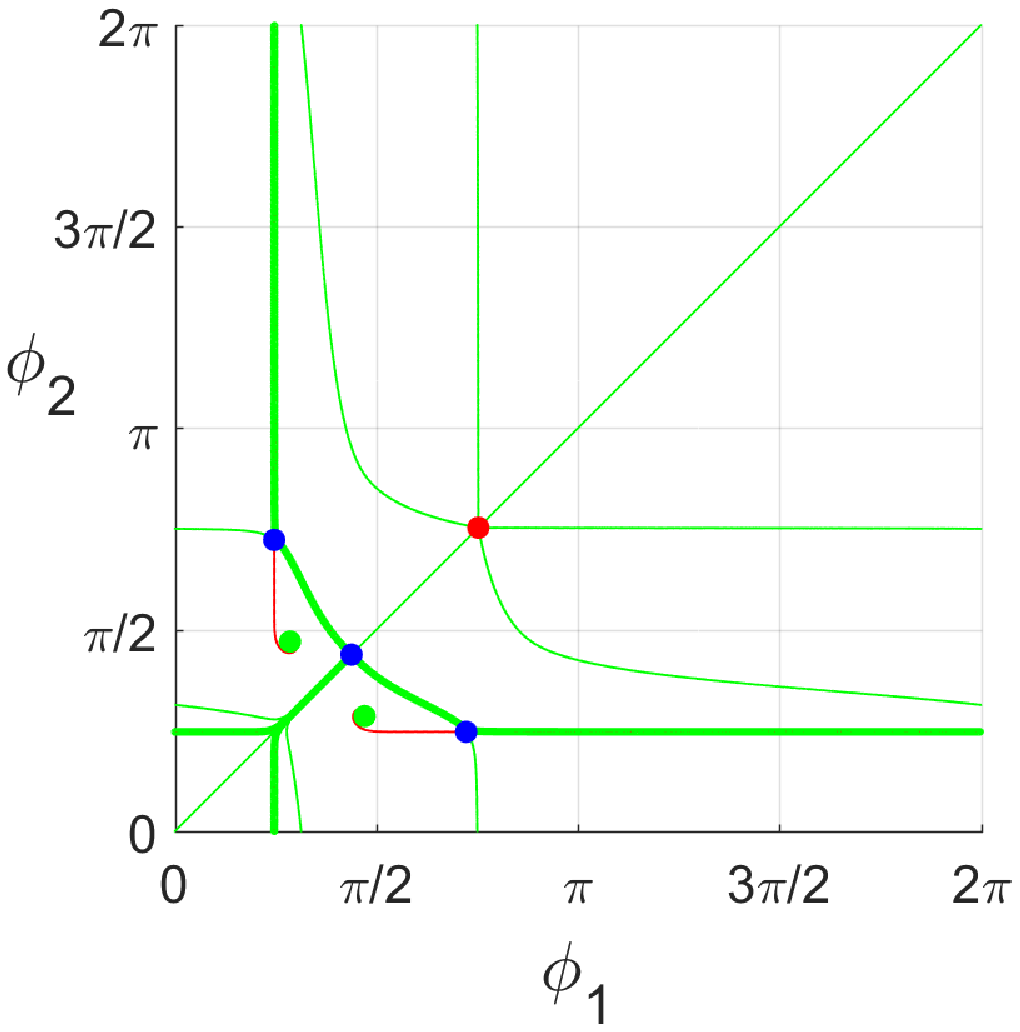}
  \caption{The fifth birth scenario of anti-phase limit cycle at the borderline between regions $C$ and $D$. $\alpha = 0.99$. (a) $\delta = 0.27$. (b) $\delta = 0.288294$. (c) $\delta = 0.288397$. (d) $\delta = 0.29$. See main text for more details. Here red/green curves correspond to unstable/stable separatrices, red dot corresponds to unstable equilibrium (unstable node), green dot -- to stable equilibrium (stable node), blue dot -- to saddle equilibrium. In (d) the green bold line corresponds to the stable anti-phase cycle.}
  \label{Born_antiphase5}
\end{figure}

The fifth scenario is also related to the emergence of the heteroclinic cycle. At its first stage, a symmetric pair of  heteroclinic trajectories appears between diagonal and non-diagonal saddles, see Fig.~\ref{Born_antiphase5}(b). Further a pair of heteroclinic trajectories between non-diagonal saddles is formed, which, along with two saddles, form the heteroclinic cycle (Fig.~\ref{Born_antiphase5}(c)). When the value of the parameter $\delta$ is further increased, this heteroclinic cycle evolves into the stable anti-phase limit cycle (Fig.~\ref{Born_antiphase5}(d)).

In the following subsection we describe how the main temporal patterns are changing in reaction to the variation of the coupling strength $d$.

\subsection{Evolution of the excitable state}
\label{sec_equilibriums_bif}

For a certain range of parameter values mentioned above the system may stay in the excitable state: self-sustained oscillations are absent, whereas the elements in the system can be excited, e.g. by an external stimulus. In the phase space in this case there is at least one stable equilibrium. We have  studied  the evolution of this basic state under the increase of coupling strength $d$. 

\begin{figure}[H]
	\centering
		(a)\includegraphics[width = 0.45\linewidth]{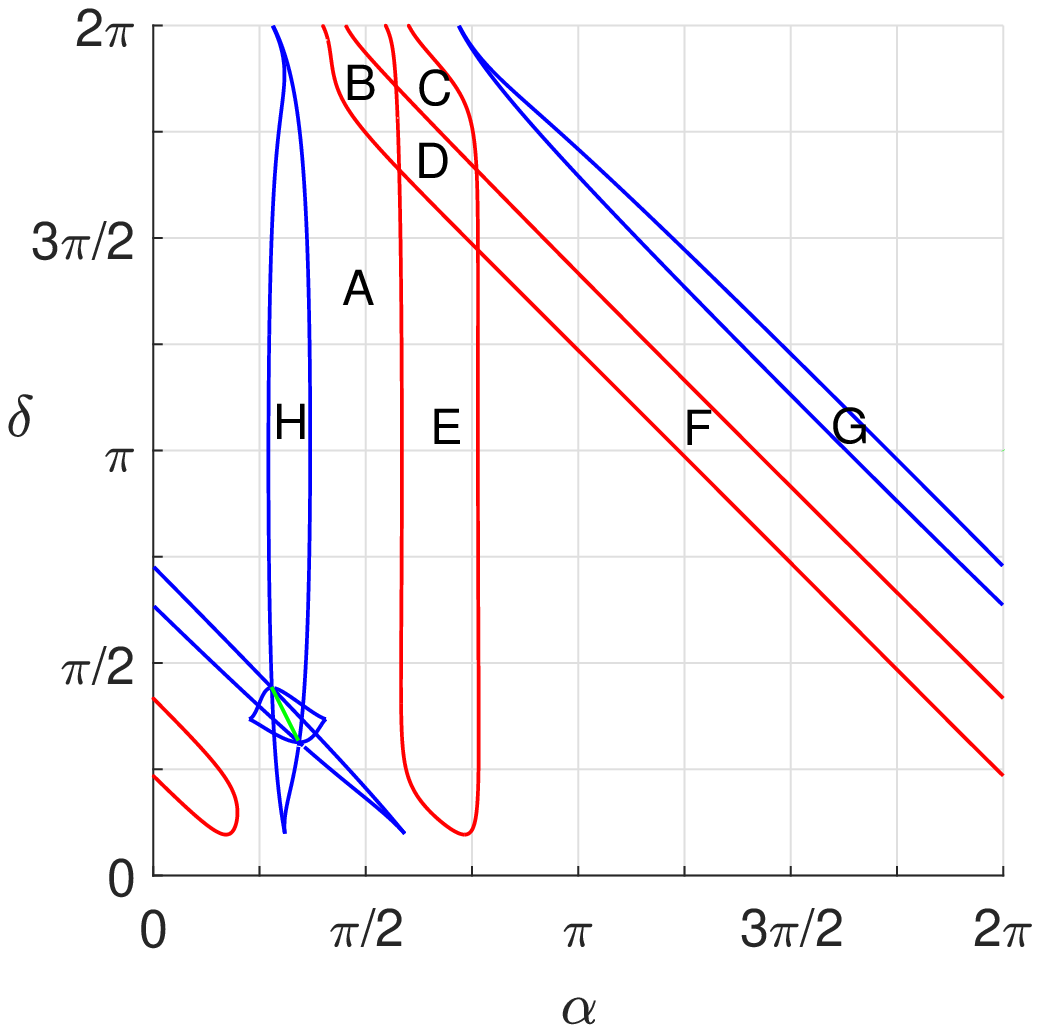}
		(b)\includegraphics[width = 0.45\linewidth]{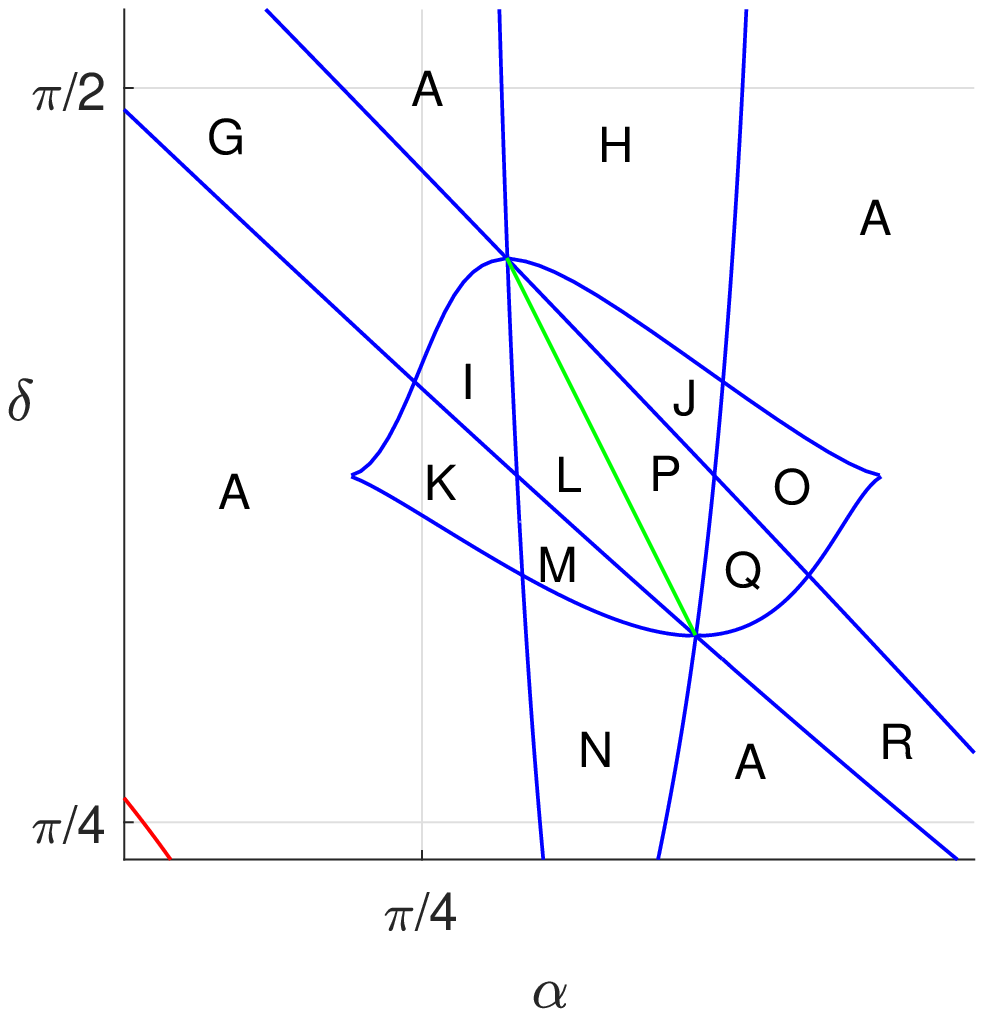}
	\caption{Bifurcation diagram on the parameter plane $(\alpha, \delta)$ and its' enlarged fragments, $d = 0.25$. All bifurcations here are bifurcations of equilibria points. Red lines correspond to pitchfork bifurcation, blue lines indicate saddle-node bifurcations, on the green line the system becomes reversible. See detailed description of regions A-Q in the text.}%, (c), (f) $d = 0.29$.}
	\label{regimes_weak_coupling}
\end{figure}

If the coupling is weak ($d < d_{th}\approx 0.3$), the structure of the phase space changes, but still only equilibria points can be observed. We divide states of rest into stable ones, unstable (they become stable in the reverse time) ones, and the saddles. Topology of torus ensures that in all configurations the number of saddle points equals the number of ``non-saddles'' (i.e., of stable and unstable states of rest). Let us list all possible setups:\\
- one stable and one unstable equilibria + two saddles (region A);\\
- two stable equilibria, one unstable equilibrium + three saddles  (regions B, F, H, N);\\
- one stable and two unstable equilibria + three saddles (regions C, E, G, R);\\
- two stable and two unstable equilibria + four saddles (region D);\\
- three stable and two unstable equilibria + five saddles (region I);\\
- two stable and three unstable equilibria + five saddles (region J);\\ 
- three stable and one unstable points + four saddles (region K);\\
- four stable and two unstable points + six saddles (region L);\\
- four stable and one unstable equilibria + five saddles (region M);\\
- one stable and three unstable equilibria + four saddles (region region O);\\
- two stable and four unstable equilibria + six saddles (region P);\\
- one stable and four unstable equilibria + five saddles (region Q).\\
This means, that, depending on the parameter values, up to four different values of the stable equilibrium membrane voltage are possible.

\subsection{Evolution of tonic spiking in dependence on the coupling strength}

Now let us list various bifurcation scenarios that, in the course of the variation of the coupling strength, lead to the onset of oscillations, including in-phase and anti-phase spiking.

\begin{figure}[H]
	\centering
		(a)\includegraphics[width = 0.3\linewidth]{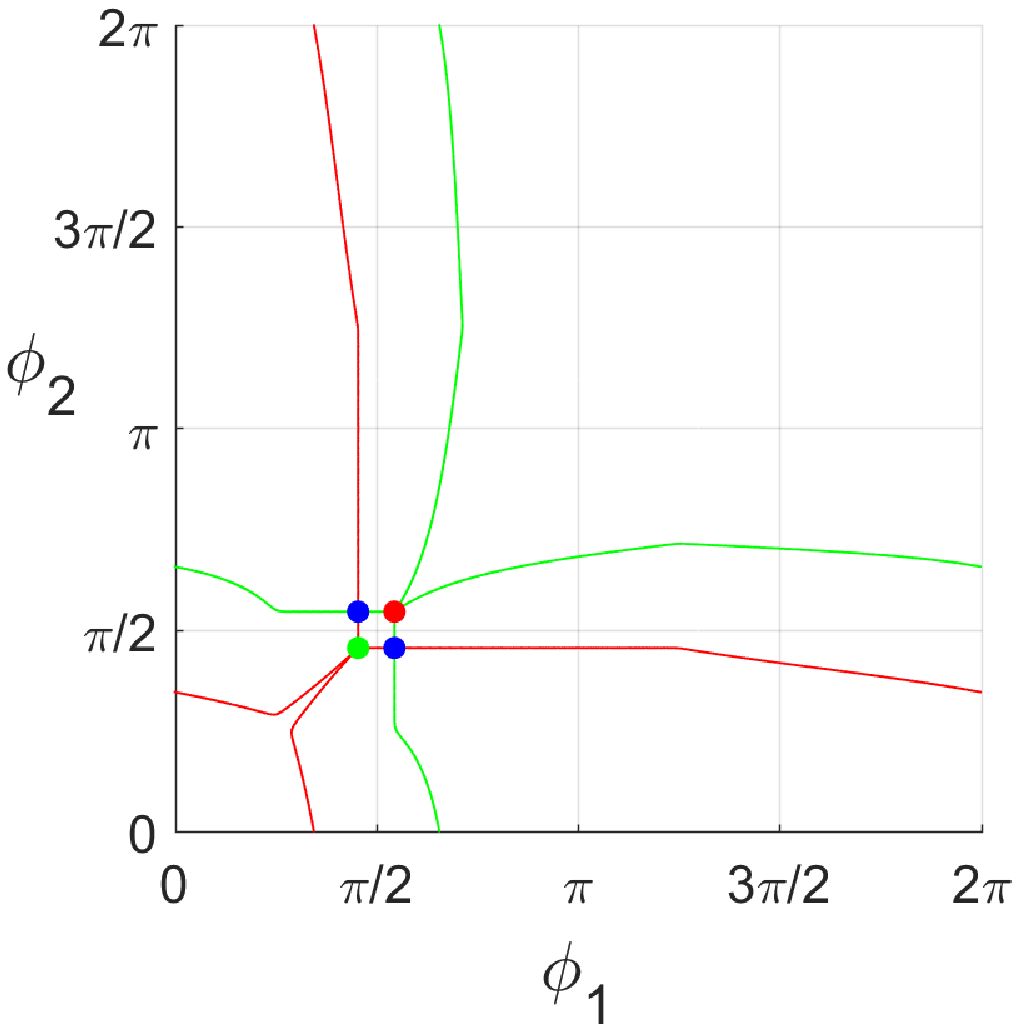}
		(b)\includegraphics[width = 0.3\linewidth]{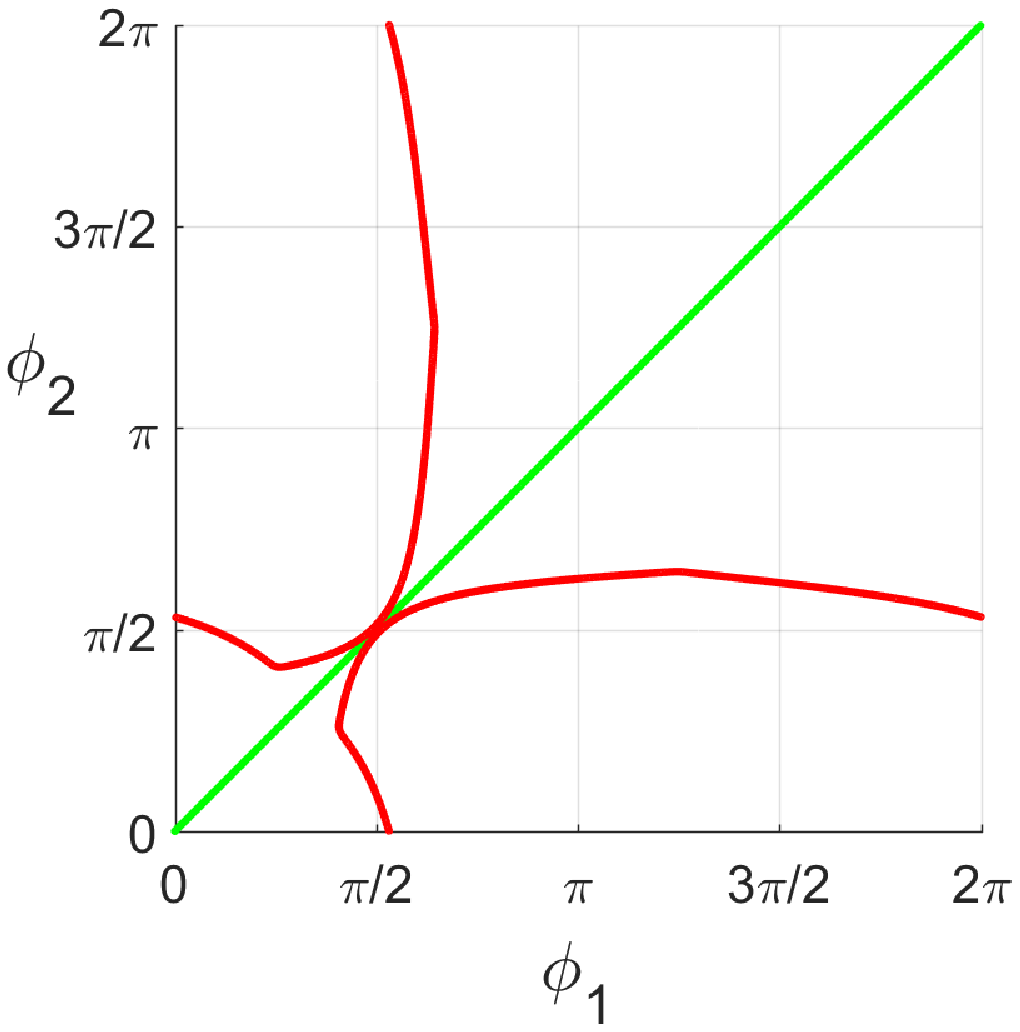}
	\caption{Birth of a stable in-phase limit cycle. Phase portraits for $\alpha =\displaystyle\frac{\pi}{4}$, $\delta = \pi$. Panel (a): $d=0.29$. Panel (b): $d=0.31$. In (a) the red/green curves show the unstable/stable separatrices. In (b) the red curve shows the unstable anti-phase limit cycle, while the green curve shows the stable in-phase limit cycle. Blue dots mark saddles, green dots show stable equilibria, red dots are unstable equilibria. See the main text for details. }
	\label{bif_anti_1}
\end{figure}

The first scenario takes place near the threshold value of the coupling strength $d_{th} = 0.3$ and is related to the onset of in-phase spiking (see Fig.~\ref{bif_anti_1}). As seen in Fig.~\ref{bif_anti_1}(a), below the threshold (here, for $d = 0.29$) two non-smooth closed invariant curves exist: the first one consists of unstable separatrices (red curves) of saddles (blue dots), saddles themselves and the stable equilibrium (green dot). This curve passes through the stable equilibrium state twice and is non-smooth at this point. The second invariant closed curve is formed by stable separatrices (green curves) of the saddles, the saddles themselves and the unstable equilibrium (red dot). When $d$ is increased, four equilibria approach each other and merge at the value $d \approx 0.3$. After the bifurcation (Fig. \ref{bif_anti_1}(b)), when the coupling strength exceeds the threshold value, e.g. for $d = 0.31$, annihilation of the equilibria is followed by formation of the in-phase stable limit cycle (green curve) and the anti-phase unstable cycle (red curve). As a result, in-phase tonic spiking is established in the system.

The condition for the birth of the anti-phase limit cycle can be approximated in the following way. 
The necessary condition for existence of limit cycles is $\gamma + d \geq 1$. Replacing the coupling function $I(\phi)$ by a piecewise constant one implies that the cycle exists if the time of motion of a phase point along the arc $(\arcsin \gamma, \pi - \arcsin \gamma)$ for the excited element does not exceed the duration of motion along the arc $(\alpha, \alpha + \delta)$ for the non-excited element:
\begin{equation}
\mathop{\int}\limits_{\arcsin \gamma}^{\pi - \arcsin \gamma} \frac{d\phi}{\gamma + d - \sin \phi} = \mathop{\int}\limits_\alpha^{\alpha + \delta} \frac{d\phi}{\gamma - \sin \phi}.
\end{equation}
This condition can be rewritten as
\begin{equation}
\begin{aligned}
\frac{2}{\sqrt{(\gamma + d)^2 - 1}}\left(\arctan{\frac{1 - (\gamma + d) \tan{\frac{\arcsin{\gamma}}{2}}}{\sqrt{(\gamma + d)^2 - 1}}} - \arctan{\frac{1 - (\gamma + d) \cot{\frac{\arcsin{\gamma}}{2}}}{\sqrt{(\gamma + d)^2 - 1}}}\right) =\\
\frac{2}{\sqrt{1 - \gamma^2}}\left(\arctanh \frac{1 - \gamma \tan{\frac{\alpha + \delta}{2}}}{\sqrt{1 - \gamma^2}} - \arctanh \frac{1 - \gamma \tan{\frac{\alpha}{2}}}{\sqrt{1 - \gamma^2}} \right)
\end{aligned}
\label{cond12}
\end{equation}

\begin{figure}[H]
	\centering
		(a)\includegraphics[width = 0.3\linewidth]{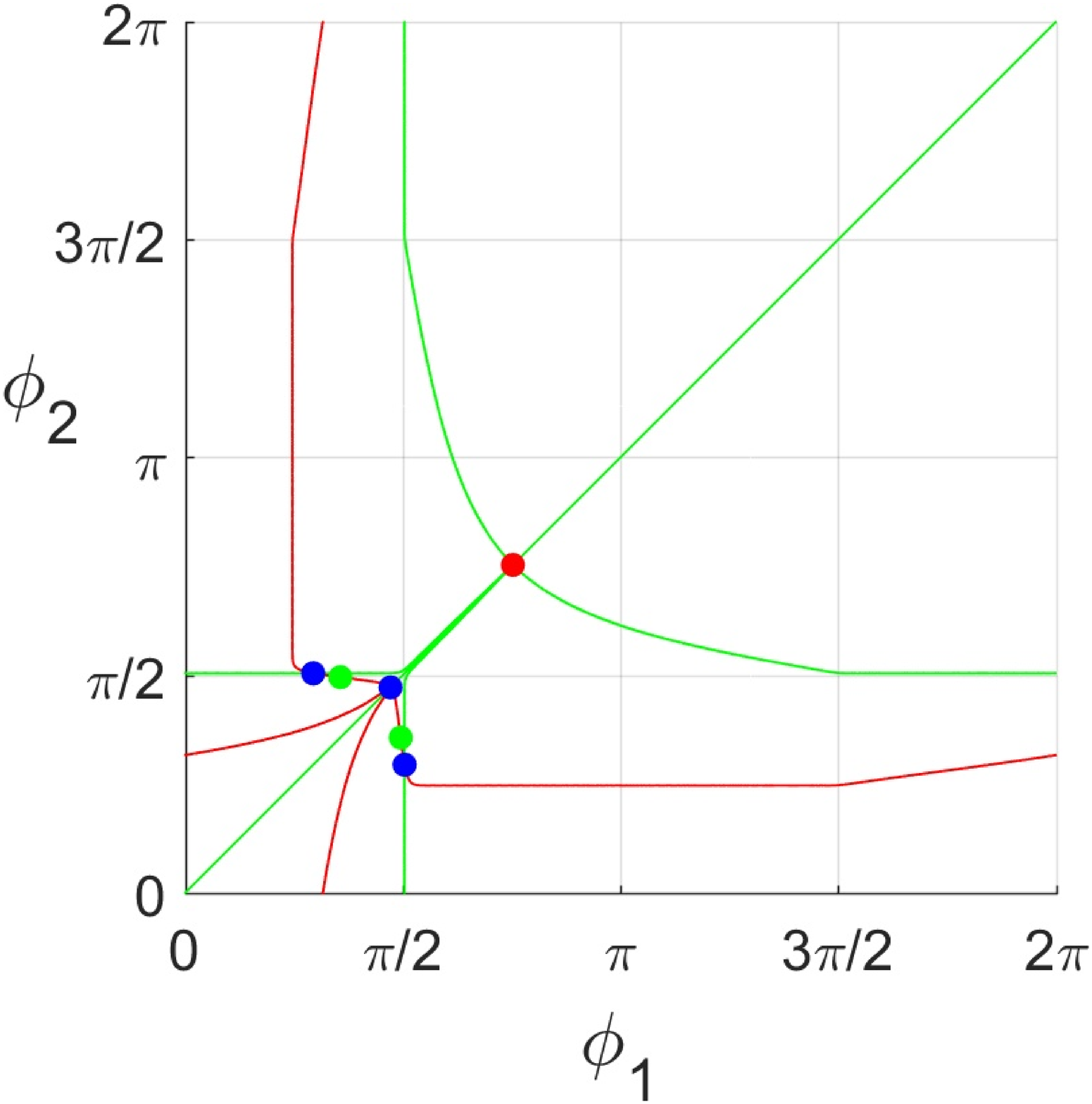}
		(b)\includegraphics[width = 0.3\linewidth]{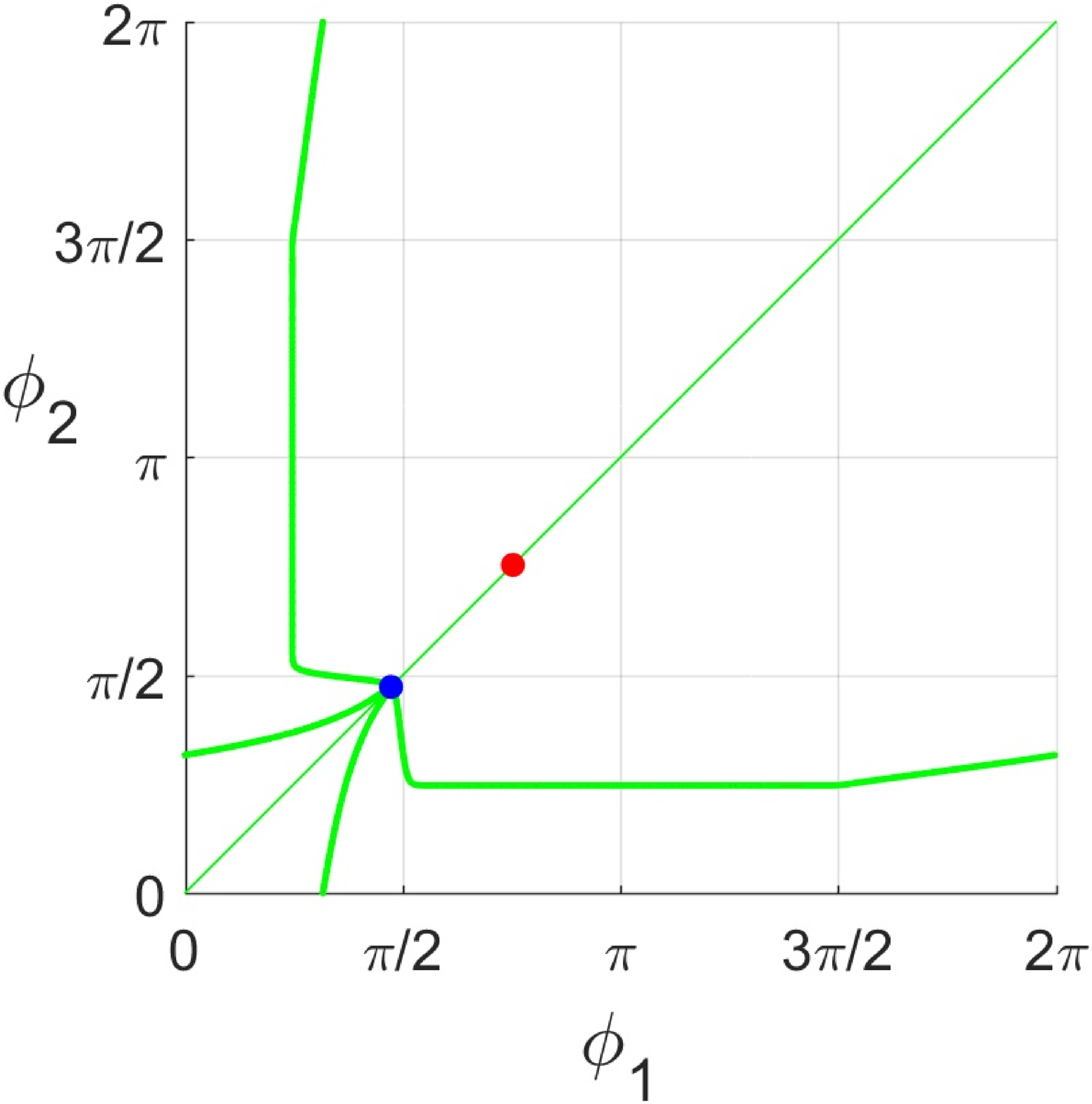}
	\caption{Birth of the stable anti-phase limit cycle. Phase portraits for $\alpha = \displaystyle\frac{3\pi}{2}$, $\delta = \pi$. (a) $d=0.2999$. (b) $d=0.301$. Red/green curves correspond to unstable/stable separatrices. In (b) the green bold curve shows the stable anti-phase limit cycle. Blue dots mark saddles, green dots --- stable equilibria, red dots --- unstable equilibria. See main text for more details. }
	\label{bif_anti_2}
\end{figure}

Bifurcation scenarios related to the appearance of anti-phase spiking pattern can be described as follows (see Fig. \ref{bif_anti_2}). 
The panel (a) shows the invariant closed curve that is formed by two non-diagonal saddles, their unstable separatrices and two stable states of rest. At $d\approx 0.3$ two saddle-node bifurcations occur on this curve, resulting in the birth of the stable anti-phase limit cycle.

\begin{figure}[H]
  \centering
  (a)\includegraphics[width = 0.33\linewidth]{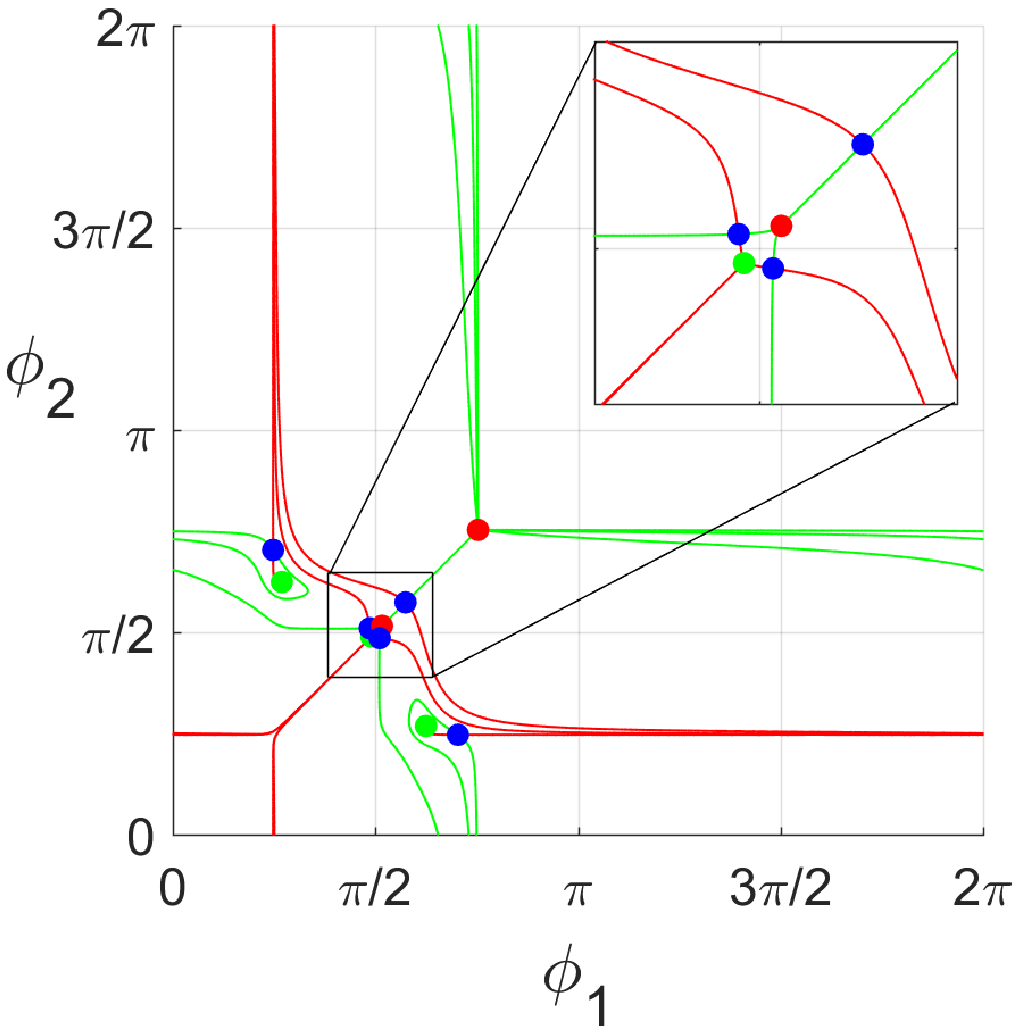}
  (b)\includegraphics[width = 0.33\linewidth]{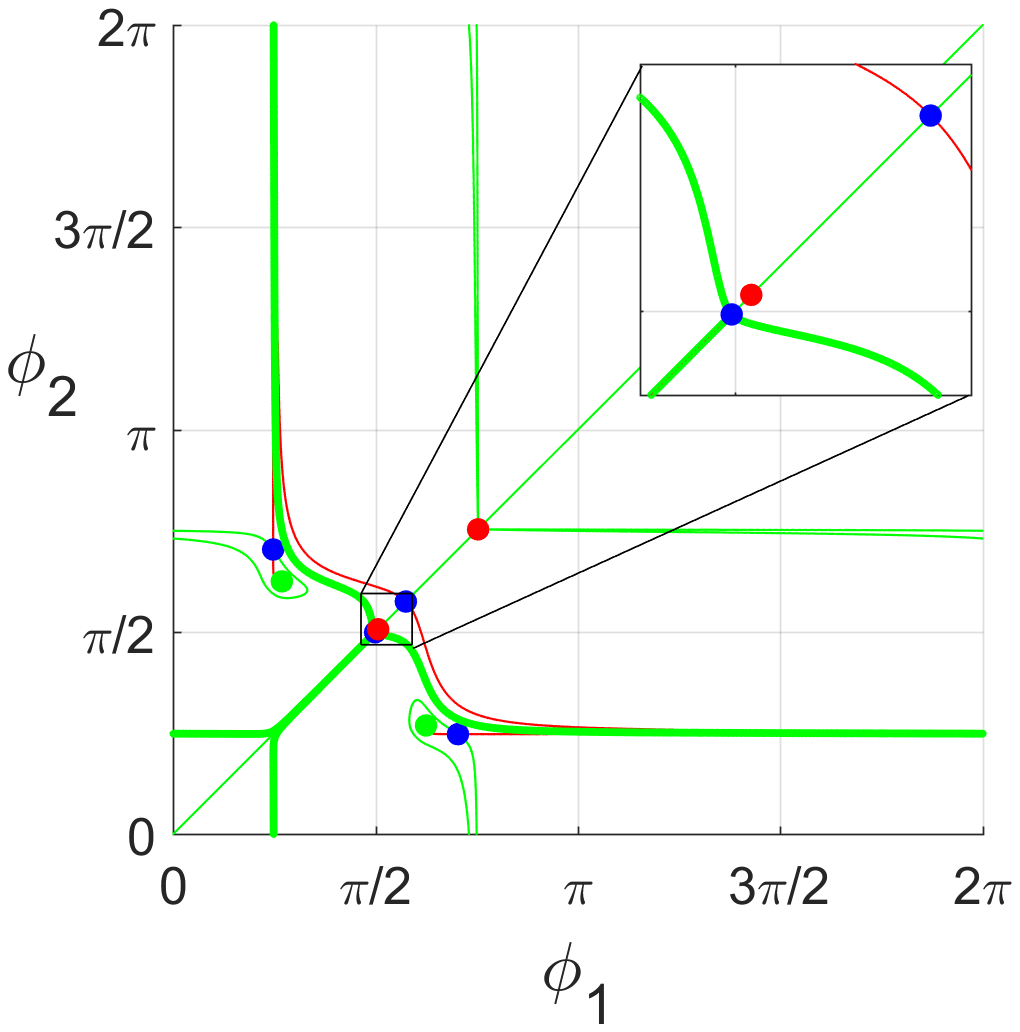}
  \caption{Birth of the stable anti-phase limit cycle for increasing coupling strength $d$ and fixed $\alpha = 0.8$, $\delta = 1.1$. Phase portraits for (a) $d = 0.2999$ (b) $d = 0.3006$. In (a) and (b) red/green curves correspond to unstable/stable separatrices. Blue dots mark saddle equilibria, green dots -- stable and red dots -- unstable ones. Bold green curve corresponds to stable anti-phase cycle. See main text for more details.}
  \label{bif_anti_d}
\end{figure}

Stages of another bifurcation scenario related to the appearance of anti-phase spiking pattern are shown in Fig.~\ref{bif_anti_d}. For coupling strength near the threshold value $d_{th}$, e.g. for $d = 0.2999$, a closed invariant curve exists. It is composed of two saddles (blue dots), their unstable separatrices (red curves), and the stable equilibrium (green dot) on the diagonal, see Fig. \ref{bif_anti_2}(a). This curve passes through the stable equilibrium twice, and is non-smooth at this point. With the increase in the value of coupling strength up to $d \approx 0.3$ the stable equilibrium undergoes a pitchfork bifurcation: it turns into a saddle which lies outside the invariant curve. Now, the invariant curve contains no equilibria and is, thereby, a limit cycle (see Fig. \ref{bif_anti_2}(b)).

\section{Conclusions}

In this study we have proposed a new phenomenological single neuron-like model and have built on its basis a model of the HCO. Constructed of two excitable neurons coupled by chemical excitatory synapses, the simple HCO model allows to conduct analytical studies; at the same time, despite its simplicity, it reflects the main properties of the biological HCO and reproduces all temporal patterns, typical for the HCO: excitable state, in-phase and anti-phase spiking. We have used bifurcation theory to obtain the mathematical description of transitions between the main types of neuron-like activity, caused by variation of the model parameters that characterize coupling. The anti-phase and in-phase spiking patterns are crucial for motor pattern generation and, according to \cite{ferrario2018bifurcations}, may be associated with swimming and synchrony patterns of spiking activity, respectively, that has been observed in a \textit{Xenopus} tadpole CPG. From the point of view of nonlinear dynamics, each of these temporal patterns corresponds to a stable periodic motion of a certain type in the phase space of the system. 

We have carried out studies of bifurcations leading to the onset of these types of neuron-like activity. On the parameter plane $(\alpha, \delta)$ where $\alpha$ characterizes the starting time of the activation of postsynaptic element and $\delta$ is responsible for duration of the couplings impact, the regions of different types of activity, such as stable in-phase and anti-phase tonic spiking, have been determined. On the parameter plane, there is also a broad region corresponding to the excitable state (quiescence), where the motif can generate activity only in response to external stimuli.

Our analysis has identified in the parameter space regions of bistability where the system demonstrates, depending on the initial conditions, both excitable and anti-phase spiking behavior, so that the same pattern generator circuit can support several types of neuron-like activity.

We have also studied transitions from excitability to spiking, caused by increase of the coupling strength $d$ from the weak coupling upwards. Obtained results, on the one hand, elucidate the origins of spiking behavior near the excitability threshold, and, on the other hand, provide deeper insights into the functions of the HCO.

%\blue{We also would like to mention that in the system under study limit cycles that are not symmetric to themselves were not observed in numerical experiments. The origin of this property is unclear. We can safely state that the reason is not the presence of symmetry in the system: one can construct very similar system with symmetry that contains pair of limit cycles that are symmetric to each other. Temporal patterns in this case are both not in-phase and not anti-phase ones.}

Remarkably, both discussed types of observed oscillatory states feature a symmetry: temporal patterns of two motif elements either coincide or are shifted with respect to each other by half of the period. Numerical search has disclosed in the parameter space neither symmetry-breaking bifurcations of the reported limit cycles nor generic periodic oscillations for which the individual states would not be related by a symmetry transformation (If present, such states are obliged to exist in pairs: the units can be interchanged). For completeness, we mention that limit cycles that are neither in-phase nor anti-phase can be encountered e.g. if $\phi$ in the denominator of the coupling function \eqref{coupling_func} is replaced by $6\phi$.

Summarizing, the proposed simple model can be used as a building block in specific complex CPG networks in a wide range of studies of motor control, dynamic memory, information processing, and decision making in animals and humans. One possible application of such studies is development of new efficient treatment of neurological diseases related to CPG arrhythmia. Another area, where these results can be helpful, concerns more efficient robot locomotion, which requires better insights in the CPG multistability \cite{kaluza2012phase}-\cite{eckert2015comparing}.

This work was partially funded by the Russian Ministry of Science and Education project \# № 14.Y26.31.0022 (studies of bifurcation scenarios) and RFBR grant \# 18-29-10068 (studies of neuronal temporal patterns).

\bibliographystyle{elsarticle-num}
\bibliography{mybibliography}

\begin{thebibliography}{10}
\expandafter\ifx\csname url\endcsname\relax
  \def\url#1{\texttt{#1}}\fi
\expandafter\ifx\csname urlprefix\endcsname\relax\def\urlprefix{URL }\fi
\expandafter\ifx\csname href\endcsname\relax
  \def\href#1#2{#2} \def\path#1{#1}\fi

\bibitem{Selverston1985}
A.~Selverston, Model neural networks and behavior, Springer Science \& Business
  Media, 2013.

\bibitem{katz2007invertebrate}
P.~S. Katz, S.~L. Hooper, Invertebrate central pattern generators, Cold Spring
  Harbor Monograph Series 49 (2007) 251.

\bibitem{mackay2002central}
M.~MacKay-Lyons, Central pattern generation of locomotion: a review of the
  evidence, Physical therapy 82~(1) (2002) 69--83.

\bibitem{guertin2013central}
P.~A. Guertin, Central pattern generator for locomotion: anatomical,
  physiological, and pathophysiological considerations, Frontiers in neurology
  3 (2013) 183.

\bibitem{matsuoka1987mechanisms}
K.~Matsuoka, Mechanisms of frequency and pattern control in the neural rhythm
  generators, Biological cybernetics 56~(5-6) (1987) 345--353.

\bibitem{pusuluri2020computational}
K.~Pusuluri, S.~Basodi, A.~Shilnikov, Computational exposition of multistable
  rhythms in 4-cell neural circuits, Communications in Nonlinear Science and
  Numerical Simulation 83 (2020) 105139.

\bibitem{selverston2000reliable}
A.~I. Selverston, M.~I. Rabinovich, H.~D. Abarbanel, R.~Elson, A.~Sz{\"u}cs,
  R.~D. Pinto, R.~Huerta, P.~Varona, Reliable circuits from irregular neurons:
  a dynamical approach to understanding central pattern generators, Journal of
  Physiology-Paris 94~(5-6) (2000) 357--374.

\bibitem{izhikevich2007dynamical}
E.~M. Izhikevich, Dynamical systems in neuroscience, MIT {P}ress, 2007.

\bibitem{sakurai2011different}
A.~Sakurai, J.~M. Newcomb, J.~L. Lillvis, P.~S. Katz, Different roles for
  homologous interneurons in species exhibiting similar rhythmic behaviors,
  Current Biology 21~(12) (2011) 1036--1043.

\bibitem{cohen1982nature}
A.~H. Cohen, P.~J. Holmes, R.~H. Rand, The nature of the coupling between
  segmental oscillators of the lamprey spinal generator for locomotion: A
  mathematical model, Journal of Mathematical Biology 13~(3) (1982) 345--369.

\bibitem{buono2001models}
P.-L. Buono, M.~Golubitsky, Models of central pattern generators for quadruped
  locomotion {I}. {P}rimary gaits, Journal of Mathematical Biology 42~(4)
  (2001) 291--326.

\bibitem{wojcik2014key}
J.~Wojcik, J.~Schwabedal, R.~Clewley, A.~L. Shilnikov, Key bifurcations of
  bursting polyrhythms in 3-cell central pattern generators, PloS one 9~(4)
  (2014) e92918.

\bibitem{jalil2013toward}
S.~Jalil, D.~Allen, J.~Youker, A.~Shilnikov, Toward robust phase-locking in
  melibe swim central pattern generator models, Chaos: An Interdisciplinary
  Journal of Nonlinear Science 23~(4) (2013) 046105.

\bibitem{hill2003half}
A.~Hill, S.~Van~Hooser, R.~Calabrese, Half-center oscillators underlying
  rhythmic movements, The handbook of brain theory and neural networks (Arbib
  M, ed) (2003) 507--510.

\bibitem{brown1911intrinsic}
T.~G. Brown, The intrinsic factors in the act of progression in the mammal,
  Proceedings of the Royal Society of London. Series B, containing papers of a
  biological character 84~(572) (1911) 308--319.

\bibitem{wang1992alternating}
X.-J. Wang, J.~Rinzel, Alternating and synchronous rhythms in reciprocally
  inhibitory model neurons, Neural computation 4~(1) (1992) 84--97.

\bibitem{terman2008reducing}
D.~Terman, S.~Ahn, X.~Wang, W.~Just, Reducing neuronal networks to discrete
  dynamics, Physica D: Nonlinear Phenomena 237~(3) (2008) 324--338.

\bibitem{alaccam2015making}
D.~Ala{\c{c}}am, A.~Shilnikov, Making a swim central pattern generator out of
  latent parabolic bursters, International Journal of Bifurcation and Chaos
  25~(07) (2015) 1540003.

\bibitem{adler1973study}
R.~Adler, A study of locking phenomena in oscillators, Proceedings of the IEEE
  61~(10) (1973) 1380--1385.

\bibitem{ermentrout1986parabolic}
G.~B. Ermentrout, N.~Kopell, Parabolic bursting in an excitable system coupled
  with a slow oscillation, SIAM Journal on Applied Mathematics 46~(2) (1986)
  233--253.

\bibitem{rubin2012explicit}
J.~E. Rubin, D.~Terman, Explicit maps to predict activation order in multiphase
  rhythms of a coupled cell network, The Journal of Mathematical Neuroscience
  2~(1) (2012) 4.

\bibitem{briggman2008multifunctional}
K.~L. Briggman, W.~Kristan~Jr, Multifunctional pattern-generating circuits,
  Annu. Rev. Neurosci. 31 (2008) 271--294.

\bibitem{destexhe1994efficient}
A.~Destexhe, Z.~F. Mainen, T.~J. Sejnowski, An efficient method for computing
  synaptic conductances based on a kinetic model of receptor binding, Neural
  computation 6~(1) (1994) 14--18.

\bibitem{korotkov2019dynamics}
A.~G. Korotkov, A.~O. Kazakov, T.~A. Levanova, G.~V. Osipov, The dynamics of
  ensemble of neuron-like elements with excitatory couplings, Communications in
  Nonlinear Science and Numerical Simulation 71 (2019) 38--49.

\bibitem{korotkov2018chaotic}
A.~G. Korotkov, A.~O. Kazakov, T.~A. Levanova, G.~V. Osipov, Chaotic regimes in
  the ensemble of fitzhhugh-nagumo elements with weak couplings,
  IFAC-PapersOnLine 51~(33) (2018) 241--245.

\bibitem{korotkov2019effects}
A.~G. Korotkov, A.~O. Kazakov, T.~A. Levanova, Effects of memristor-based
  coupling in the ensemble of fitzhugh--nagumo elements, The European Physical
  Journal Special Topics 228~(10) (2019) 2325--2337.

\bibitem{ferrario2018bifurcations}
A.~Ferrario, R.~Merrison-Hort, S.~R. Soffe, W.-C. Li, R.~Borisyuk, Bifurcations
  of limit cycles in a reduced model of the xenopus tadpole central pattern
  generator, The Journal of Mathematical Neuroscience 8~(1) (2018) 10.

\bibitem{kaluza2012phase}
P.~Kaluza, T.~Cioac{\u{a}}, Phase oscillator neural network as artificial
  central pattern generator for robots, Neurocomputing 97 (2012) 115--124.

\bibitem{eckert2015comparing}
P.~Eckert, A.~Spr{\"o}witz, H.~Witte, A.~J. Ijspeert, Comparing the effect of
  different spine and leg designs for a small bounding quadruped robot, in:
  2015 IEEE International Conference on Robotics and Automation (ICRA), IEEE,
  2015, pp. 3128--3133.

\end{thebibliography}

\end{document}